\documentclass[reqno,12pt]{amsart}

\usepackage{amssymb,amsmath,amsthm,stmaryrd}
\usepackage{pstricks,pst-node,multido,pst-tree,pstcol}
\usepackage{euscript}
\usepackage{calc,ifthen}
\usepackage{array,hhline}
\usepackage{paralist}
\usepackage{algorithmic,algorithm}
\usepackage{float,caption}
\usepackage{flafter}
\usepackage{threeparttable}
\usepackage{url}
\usepackage{enumerate}
\usepackage{graphics,color}
\DeclareMathAlphabet{\mathcal}{OMS}{ztmcm}{m}{n}

\theoremstyle{plain}
\newtheorem{theorem}{Theorem}[section]
\newtheorem{proposition}[theorem]{Proposition}
\newtheorem{lemma}[theorem]{Lemma}
\newtheorem{corollary}[theorem]{Corollary}

\theoremstyle{definition}
\newtheorem{definition}[theorem]{Definition}
\newtheorem{example}[theorem]{Example}
\newtheorem{remark}[theorem]{Remark}

\newcommand{\vc}[1]{\mathord{\mbox{$\vcenter{\hbox{#1}}$}}}

\newcommand\refcond[1]{($\mathrm{R}_{\ref{#1}}$)}
\newcommand{\T}{\EuScript{T}}

\newcommand{\PT}{\EuScript{PT}}

\newcommand{\Z}{{\mathbb Z}}
\newcommand{\N}{\mathbb{N}}

\newcommand{\FDD}{DD}
\newcommand{\FDT}{DT}

\newcommand{\field}{\Bbbk}
\newcommand{\BT}{\operatorname{\mathsf{BT}}}
\newcommand{\DT}{\operatorname{\mathcal{D\!T}}}
\newcommand{\DD}{\operatorname{\mathcal{DD}}}

\DeclareMathOperator{\degnode}{deg_{\mathsf{node}}}
\DeclareMathOperator{\degangle}{deg_{\mathsf{angle}}}

\DeclareMathOperator{\assocprod}{\diamond}
\DeclareMathOperator{\assocprodlambda}{\diamond_\lambda}
\DeclareMathOperator{\astlambda}{\ast_\lambda}

\newcounter{xx}
\newcommand{\ver}{\ifthenelse{\value{xx} = 0}{\vrule}{}}

\newcommand{\xyinc}{\ar@{^{(}->}}

\newcommand{\onto}{\twoheadrightarrow} 
\newcommand{\map}[1]{\xrightarrow{#1}}
\newcommand{\inc}{\hookrightarrow}

\newcommand{\tii}{\T_{\infty,\infty}}
\newcommand{\tiw}{\T_{\infty,2}}
\newcommand{\twi}{\T_{2,\infty}}
\newcommand{\tww}{\T_{2,2}}
\newcommand{\txy}{\T_{i,j}}

\newcommand\BinT{\EuScript{BT}} % Binary trees
\newcommand\SP{\EuScript{SP}} % Schr\"oder paths
\newcommand\SPz{\EuScript{SP}^0} % Schr\"oder paths, D in diagonal
\newcommand\SPp{\EuScript{SP}^+} % Schr\"oder paths, D not in diagonal
\newcommand\RP{\EuScript{RP}} % 'Rota-Baxter' paths
\newcommand\RPz{\EuScript{RP}^0} % 'Rota-Baxter' paths
\newcommand\RPp{\EuScript{RP}^+} % 'Rota-Baxter' paths
\newcommand\CP{\EuScript{CP}} % Catalan paths
\newcommand\MP{\EuScript{MP}} % Motzkin paths
\newcommand\RMP{\EuScript{RMP}} % 'Rota-Baxter' Motzkin paths
\newcommand\hMP{\mathord{\MP_h}} % Stanley colored MP
\newcommand\huMP{\mathord{\MP_{hu}}} % up-horiz. colored MP

\newcommand{\betaii}{\beta_{\infty,\infty}}
\newcommand{\betaiw}{\beta_{\infty,2}}

\newcommand{\betaww}{\beta_{2,2}}
\newcommand{\betaxy}{\beta_{i,j}}

\newcommand{\f}{\mathsf{b}}
\newcommand{\fxy}{\f_{i,j}}
\newcommand{\fii}{\f_{\infty,\infty}}
\newcommand{\fiw}{\f_{\infty,2}}
\newcommand{\fwi}{\f_{2,\infty}}
\newcommand{\fww}{\f_{2,2}}

\newcommand{\dt}{\mathsf{dt}}

\newcommand{\Gii}{\mathsf{B}_{\infty,\infty}}
\newcommand{\Giw}{\mathsf{B}_{\infty,2}}
\newcommand{\Gwi}{\mathsf{B}_{2,\infty}}
\newcommand{\Gww}{\mathsf{B}_{2,2}}

\newcommand{\R}{\mathcal{B}}
\newcommand{\Rii}{\R_{\infty,\infty}}
\newcommand{\Rxy}{\R_{i,j}}
\newcommand{\Riw}{\R_{\infty,2}}
\newcommand{\Rwi}{\R_{2,\infty}}
\newcommand{\Rww}{\R_{2,2}}

\newcommand{\Mori}{\EuScript{M}_{\infty}}
\newcommand{\Moriw}{\EuScript{M}_{2}}
\newcommand{\Mi}{M_\infty}
\newcommand{\Mw}{M_2}
\newcommand{\betai}{\beta_\infty}
\newcommand{\betaw}{\beta_2}
\newcommand{\pii}{\pi_\infty}
\newcommand{\piw}{\pi_2}

\newcommand{\mi}{\mathsf{m}_{\infty}}
\newcommand{\mw}{\mathsf{m}_{2}}

\newcommand{\F}{B}
\newcommand{\Fii}{\F_{\infty,\infty}}
\newcommand{\Fiw}{\F_{\infty,2}}
\newcommand{\Fwi}{\F_{2,\infty}}
\newcommand{\Fww}{\F_{2,2}}
\newcommand{\Fxy}{\F_{i,j}}

\newcommand{\Gr}{\mathsf{G}}
\newcommand{\Grxy}{\mathsf{G}_{i,j}}
\newcommand{\Grii}{\mathsf{G}_{\infty,\infty}}

\renewcommand\H{\mathsf{H}}
\newcommand{\Hxy}{\mathsf{H}_{i,j}}
\newcommand{\Hii}{\mathsf{H}_{\infty,\infty}}

\newcommand\treetopath[1]{\operatorname{\mathsf{TreeToPath}}(#1)}
 
\SpecialCoor
\newlength\levelsep 
\newlength\treesepone 
\newlength\treeseptwo 
\newlength\lblsep 
\newlength\lblheight 
\newlength\temp

\psset{linewidth=0.5pt,
       dotsize=1.8pt 1,
       unit=1em}

\def\angleedge#1#2#3{\ncline{c-c}{#2}{#3}\ncput[npos=0.7]{\pnode{#1}}}
\def\angleedgeX#1#2#3{\ncline{c-c}{#2}{#3}\ncput[npos=0.8]{\pnode{#1}}}

\newcommand\children[5]{% name, {,options}, howmany, dist, levelsep
  \pstree[levelsep=#5,treesep=#4#2]{\Tdot\pnode{#1Root}}%
  {\multido{\ia=1+1}{#3}{{\Tdot[edge=\angleedge{I#1Child\ia}]\pnode{#1Child\ia}}}}%
}

\newcommand\childrenX[5]{% name, {,options}, howmany, dist, levelsep
  \pstree[levelsep=#5,treesep=#4#2]{\Tdot\pnode{#1Root}}%
  {\multido{\ia=1+1}{#3}{{\Tdot[edge=\angleedgeX{I#1Child\ia}]\pnode{#1Child\ia}}}}%
}

\newcommand\tree[1]{%
  \mathord{\hborder{1.5pt}{\vertical{#1}}}%
}

\newcommand\angledX[3]{%
  \ncline[linestyle=none]%
  {-}{I#1}{I#2}\ncput[labelsep=0pt]%
    {\ensuremath{\scriptscriptstyle{\textcolor{blue}{#3}}}}}

\newcommand\angledG[3]{%
  \ncarc[linewidth=0.3pt,linecolor=blue,arcangle=-20,nodesep=5pt]%
  {-}{I#1}{I#2}\naput[labelsep=1pt]%
    {\ensuremath{\scriptscriptstyle{\textcolor{blue}{#3}}}}}

\newcommand{\angled}[3]{%
  \ncarc[linewidth=0.3pt,linecolor=blue,arcangle=-50,nodesep=1.5pt]%
  {-}{I#1}{I#2}\naput[labelsep=1pt]%
    {\ensuremath{\scriptscriptstyle{\textcolor{blue}{#3}}}}}

\newcommand\lbl[5][c]{% ref, node, vnodesep, hnodesep, stuff
  \rput[#1](#2){\vborder{#3}{\hborder{#4}{\ensuremath{\scriptstyle{\textcolor{red}{#5}}}}}}}

\newcommand\point{\mathord{\vertical{$\;\psdots[dotsize=2pt 1](0,0)\;$}}}

\newcommand{\Y}[2]{%
  \preY{}{#1}{#2}%
}

\newcommand{\preY}[3]{% {,options}, root label, angle label
  \settoheight\temp{$\scriptstyle{#2}$}
  \tree{\children{A}{#1,xbbh=.5\temp}{2}{12pt}{2.5ex}%
      \lbl[l]{ARoot}{0pt}{\lblsep}{#2}%
      \angled{AChild1}{AChild2}{#3}%
  }
}

\newcommand{\preYY}[3]{% {,options}, root label, angle label
  \settoheight\temp{$\scriptstyle{#2}$}
  \tree{\children{A}{#1,xbbh=.5\temp}{2}{10pt}{2.5ex}%
      \lbl[l]{ARoot}{0pt}{\lblsep}{#2}%
  }
}

\newcommand{\generator}[1][1]{\Y{0}{#1}}
\newcommand{\betagenerator}[1][1]{\Y{1}{#1}}
\newcommand{\gendend}{\preYY{}{}{}}

\newenvironment{descr}[1][55pt]
  {\begin{list}{}{ %
        \setlength\labelwidth{#1}%
        \setlength\labelsep{0pt}%
        \setlength\itemsep{0.5ex}%
        \setlength\leftmargin{#1}}}
  {\end{list}}

\newlength\verticalwidth
\newsavebox\verticalbox
\newcommand\vertical[1]{%
  \settowidth\verticalwidth{#1}%
  \parbox[c]{\verticalwidth}{#1}}

\newcommand\vborder[2]{%
  #2\raisebox{#1}{\vphantom{#2}}\raisebox{-#1}{\vphantom{#2}}}
\newcommand\hborder[2]{%
  \makebox[\width+(#1)*2][c]{#2}}

\newcommand{\rb}{Baxter}
\newcommand{\rba}{\rb\ algebra}

\title{Combinatorics of the free Baxter algebra}
\author[M. Aguiar and W. Moreira]{Marcelo Aguiar and Walter Moreira}
\thanks{We
    thank Kurusch Ebrahimi-Fard for an explanation of the
    paper~\cite{EG}, which led us to the results of this paper.}

\address{
Department of Mathematics\\
Texas A\&M University\\
College Station, TX 77843, USA}
\email{maguiar@math.tamu.edu}
\urladdr{\url{http://www.math.tamu.edu/~maguiar}}

\address{Department of Mathematics\\ 
Texas A\&M University\\
College Station, TX 77843, USA}
\email{wmoreira@math.tamu.edu}
\urladdr{\url{http://www.math.tamu.edu/~wmoreira}}

\date{March 16, 2007}
\begin{document}

\begin{abstract}
We study the free (associative, non-commutative) \rba\ on one generator. The first explicit description of this object is due to Ebrahimi-Fard and Guo. We provide an alternative description in terms of a certain class of trees, which form a linear basis for this algebra. We use this to treat other related cases, particularly that in which the \rb\ map is required to be quasi-idempotent, in a unified manner. Each case corresponds to a different class of trees.

Our main focus is on the underlying combinatorics. In several cases, we  provide bijections between our various classes of trees 
and more familiar combinatorial objects including certain Schr\"oder paths and Motzkin paths.
We calculate the dimensions of the homogeneous components of these algebras (with respect to a bidegree related to the number of nodes and the number of angles in the trees) and the corresponding generating series. An important feature is that the
combinatorics is captured by the idempotent case; the others are obtained from this case by various binomial transforms.
We also relate free \rb\ algebras to Loday's dendriform trialgebras and dialgebras. We show that the free dendriform trialgebra (respectively, dialgebra) on one generator
embeds in the free \rba\ with a quasi-idempotent map (respectively, with a quasi-idempotent map and an idempotent generator). This refines results of
Ebrahimi-Fard and Guo.
\end{abstract}

\maketitle

\reversemarginpar

%\psset{radius=1.5pt,%
%  levelsep=4ex,%
%  treesep=9pt,%
%  labelsep=3pt,
%  showbbox=false,%
%%  showbbox=true,%
%  linewidth=0.5pt,%
%  unit=1em
%}

\section{Introduction}

A {\em \rba} (also called {\em Rota-Baxter algebra} in some of the recent literature) is a pair $(A,\beta)$ consisting of an associative algebra $A$ and a linear map $\beta:A\to A$ satisfying
\[ \beta(a)\beta(b) = \beta\bigl(\beta(a)b + a\beta(b) + \lambda ab\bigr),\]
where $\lambda$ is a fixed scalar.
Interest in these objects originated in work of Baxter~\cite{B}. 
Constructing the free \rba\ in explicit terms amounts to describing all consequences of the above identity.  

Rota gave the first description of the free {\em commutative} \rba\ ~\cite{R69}, by providing an embedding into an explicit \rba. Cartier then obtained
an intrinsic description~\cite{Ca}. For other references to early work,
see~\cite{R1,R2}.
More recently, Guo and Keigher described the adjoint functor to the forgetful functor from the category of commutative  \rba s  to the category of commutative algebras~\cite{GK}. 

It is natural to consider the possibly more challenging task of
constructing the free \rba, not necessarily commutative. In recent
interesting work, Ebrahimi-Fard and Guo have successfully tackled this
problem~\cite{EG}; they have in fact constructed the adjoint functor
to the forgetful functor from the category of (associative) \rba s to
the category of (associative) algebras. As it turns out, there is not
much loss of generality in concentrating in the case of one generator
$x$, which we do from now on.  The construction in~\cite{EG} involves
a certain class of words on the symbols $x$ and $\beta(x)$. This
choice of combinatorial structure makes the description of the
algebraic structure rather involved and lengthy.

In this paper we  provide a  simpler
description of this algebra, by making use of a different combinatorial structure (decorated trees) and of an appropriate notion of grafting for these objects. 
 We have learned that the authors of~\cite{EG}
 were aware of this possibility, and plan to present their results in~\cite{EG2}. Another paper in preparation with related results to ours is~\cite{GS}.

The use of decorated trees makes our construction very reminiscent of the constructions of the free dendriform dialgebra of Loday~\cite{L} 
and  of the free dendriform trialgebra 
 of Loday and Ronco~\cite{LR}. In addition, it
 allows us to present a unified construction of the free \rba\ and of three closely related algebras; namely, that in which the generator $x$ is assumed
to be idempotent $(x^2=x)$,  that in which the map $\beta$ is assumed to be
quasi-idempotent $(\beta^2=-\lambda\beta)$, and that in which both assumptions are made. 
 We refer to any of these as a free \rba\ (of the appropriate kind) and
denote them by $\Fxy^\lambda$, where the subindices $i,j\in\{2,\infty\}$ distinguish between the various cases. They are related by a commutative diagram of surjective morphisms of \rba s as follows:
\[\vertical{\ensuremath{\psset{nodesep=1pt}%
      \begin{psmatrix}[rowsep=8pt,colsep=10pt]
        & \Fii^\lambda \\
        \Fiw^\lambda & & \Fwi^\lambda \\
        & \Fww^\lambda
        \ncline{->>}{1,2}{2,3} \ncline{->>}{1,2}{2,1}
        \ncline{->>}{2,1}{3,2} \ncline{->>}{2,3}{3,2}
      \end{psmatrix}%
    }}%
\]
\enlargethispage{5pt}
The free \rba s $\Fwi^{\lambda}$ (in which the generator is assumed to be idempotent) and  $\Fii^{\lambda}$ (in which no assumptions are made)
are covered by the adjoint construction of~\cite{EG}. The algebras $\Fww^\lambda$ and $\Fiw^\lambda$ (in which the map is assumed to be quasi-idempotent) %are not treated in those references and they
 constitute the main focus of our work. For our purposes these cases appear to be more fundamental, as explained in the next three paragraphs.

One of our goals is to calculate the dimensions of the homogeneous components of the algebras $\Fxy^\lambda$, and the corresponding generating series. An important feature is that the combinatorics is captured by the idempotent case:
the generating series for the algebras $\Fxy^\lambda$ are {\em binomial transforms} of the generating series for the algebra $\Fww^\lambda$. 
We provide explicit formulas for the dimensions  of the homogeneous components of the algebras $\Fww^\lambda$ and $\Fiw^\lambda$, and on the way to these results we provide several bijections between the classes of decorated trees that form linear bases of these algebras and more familiar combinatorial objects, such as planar rooted trees, Schr\"oder paths, and Motzkin paths. For a summary of the most important combinatorial results,
see Table~\ref{T:dimensions}.

Another goal is to clarify the connections between free \rba s and free dendriform dialgebras and trialgebras. Dendriform dialgebras and trialgebras were introduced by Loday~\cite{L} and
Loday and Ronco~\cite{LR}. A connection between these objects and
\rba s was observed in~\cite{poisson,E}: any  \rba\ with $\lambda=1$ can be turned into a dendriform trialgebra and any \rba\ with $\lambda=0$ can be turned into a dendriform dialgebra. This gives rise to morphisms of \rba s from the free dendriform trialgebra on one generator to $\Fii^1$ and from the free dendriform dialgebra on one generator to $\Fii^0$. Ebrahimi-Fard and Guo showed that these maps are injective~\cite{EG}. We show here that in fact the  free dendriform trialgebra embeds in $\Fiw^1$ and the free dendriform dialgebra embeds
in $\Fww^0$.

We also discuss algebras $A$ equipped with an idempotent endomorphism of algebras $\beta$. Such a pair $(A,\beta)$ is a \rba\ with $\lambda=-1$, so choosing an element of $A$ determines  a morphism $\Fiw^{-1}\to A$ of \rba s.
We construct the free object on one generator in this category and describe the canonical morphism from $\Fiw^{-1}$ in explicit terms. We also provide the analogous results for the case of idempotent generators.

Decorated trees are introduced in Section~\ref{S:decorated}, and the notion of grafting, which is central for the construction of the free \rba s, is discussed in~\ref{S:grafting}. The construction is carried out in Section~\ref{S:freeRB}, where
we provide a complete concise proof of the universal property of the algebras $\Fxy^\lambda$ (Proposition~\ref{mainthm}). Section~\ref{S:combinatorics} contains the combinatorial results; though our motivation is algebraic, these results are interesting on their own, 
and they can be read separately from the rest.
Section~\ref{S:comb-trees} presents various kinds of combinatorial objects and
then puts them in bijection with the linear bases of the free \rba s. These results are used to calculate the dimensions of the homogeneous components of the free \rba s  in Section~\ref{S:comb-dim}, as well as  the generating series in~\ref{S:generating}.  
Algebras with an idempotent endomorphism and their connection to \rba s are
discussed in Section~\ref{S:idem}. The connection with dendriform trialgebras and dialgebras and the embedding results are given in Section~\ref{S:dendri}.
The appendix contains two algorithms used to set up some of the bijections of
Section~\ref{S:combinatorics}.

\subsection*{Notation} We work over a commutative ring $\field$. By vector space we mean free $\field$-module. All spaces and algebras are over $\field$.
All algebras are associative, but not necessarily unital.

The set $\Z^+$ is the set of positive integers and $\N=\Z^+\cup\{0\}$.

\section{Free \rba s on one generator}
\label{sectionRotaBaxter}

Let $A$ be an algebra, $\lambda\in\field$, and $\beta:A\to A$ 
a linear map satisfying
\begin{equation}\label{E:rb}
  \beta(a)\beta(b) = \beta\bigl(\beta(a)b + a\beta(b) + \lambda ab\bigr)
\end{equation}
for all $a,b\in A$. The map $\beta$ is called a {\em \rb\ operator} and
the pair $(A,\beta)$ is called a {\em \rba\ of weight $\lambda$}.
In this case, defining
\begin{equation}\label{E:def-ast}
a\astlambda b=\beta(a)b + a\beta(b) + \lambda ab
\end{equation}
one obtains a new associative operation on $A$.

The free \rba\  was constructed by Ebrahimi-Fard and 
Guo~\cite{EG}. Below we provide another
description of the free \rba\ on one generator, as well as of three related algebras in which either the generator $x$ is assumed to be idempotent:
\begin{equation}\label{E:gen-idem}
x^2=x,
\end{equation}
or the \rb\ map $\beta$ is assumed to be quasi-idempotent:
\begin{equation}\label{E:map-idem}
\beta^2=-\lambda\beta.
\end{equation}
 Our description is in terms of decorated trees,
as discussed in Section~\ref{S:decorated} below. This allows us to provide
simpler definitions of the product in these algebras and of the \rb\ maps.
It also proves useful in calculating the dimensions of the
homogeneous components of these algebras, see Section~\ref{S:comb-dim}.

\begin{remark}\label{R:exceptional}
One may wonder about imposing the relation
\[\beta^2=\mu\beta\]
where $\mu\in\field$ is some scalar other than $-\lambda$. 
In this case, additional relations follow from~\eqref{E:rb} and the above relation, such as $\beta\bigl(a\beta(b)\bigr)=\beta\bigl(\beta(a)b\bigr)=0$ and
$\beta(a)\beta(b)=\lambda\beta(ab)$. 
This leads to
three different constructions (according to whether $\lambda=0$ or $\mu=0$)
which we do not treat in this paper. 
\end{remark}

\subsection{Decorated trees}\label{S:decorated}

We define the sets that are
going to be bases as vector spaces of the free \rba s.

Consider a rooted planar tree $t$. A node of $t$ is a  {\em leaf} if it has no children, otherwise it is an {\em internal node}.
 An {\em angle} of $t$ is the
sector between two consecutive children of an internal node. We
decorate $t$ by writing positive integers in the angles and
non-negative integers on the internal nodes:
\begin{equation*}
\text{rooted planar tree}\quad
\vertical{\children{A}{,xbbd=\levelsep}{3}{\treesepone}{\levelsep}%
  \rput[t](AChild2){%
    \children{B}{}{2}{\treesepone}{\levelsep}}}
\enspace,
\qquad%
%%%
\text{decorated rooted planar tree}\quad%
\vertical{\children{A}{,xbbd=\levelsep}{3}{\treesepone}{\levelsep}%
  \rput[t](AChild2){%
    \children{B}{}{2}{\treesepone}{\levelsep}}%
  \lbl[bl]{ARoot}{0pt}{\lblsep}{1}%
  \lbl[l]{BRoot}{0pt}{\lblsep}{\raisebox{-7pt}{\ensuremath{\scriptstyle 4}}}%
  \angled{AChild1}{AChild2}{2}%
  \angled{AChild2}{AChild3}{1}%
  \angled{BChild1}{BChild2}{5}%
}
\enspace.
\end{equation*}
Let $\tii$ be the set consisting of all decorated rooted planar trees
satisfying the following conditions:
\begin{enumerate}[($\mathrm{R}_1$)]
\item \label{InternalFertility} Every internal node has at least two children.
\item \label{IntermediateChildren}Among the children of each node, only the leftmost and rightmost
  children can be leaves.
  \item \label{RootLabel} Only the root may have label
  $0$; all other internal nodes must be labeled with positive integers.
\end{enumerate}
For example,  among the following trees,
\begin{equation*}
t_1:\quad
\vertical{\children{A}{}{3}{\treesepone}{\levelsep}%
  \lbl[bl]{ARoot}{0pt}{\lblsep}{2}%
  \angled{AChild1}{AChild2}{2}%
  \angled{AChild2}{AChild3}{1}}
\quad\qquad
%%%
t_2:\quad
\vertical{\children{A}{,xbbd=\levelsep}{3}{\treesepone}{\levelsep}%
  \rput[t](AChild2){%
    \children{B}{}{2}{\treesepone}{\levelsep}}%
  \lbl[bl]{ARoot}{0pt}{\lblsep}{1}%
  \lbl[l]{BRoot}{0pt}{\lblsep}{\raisebox{-7pt}{\ensuremath{\scriptstyle 0}}}%
  \angled{AChild1}{AChild2}{2}%
  \angled{AChild2}{AChild3}{1}%
  \angled{BChild1}{BChild2}{2}%
}
\quad\qquad
%%%
t_3:\quad
\vertical{\children{A}{,xbbd=\levelsep}{3}{\treesepone}{\levelsep}%
  \rput[t](AChild2){%
    \children{B}{}{2}{\treesepone}{\levelsep}}%
  \lbl[bl]{ARoot}{0pt}{\lblsep}{1}%
  \lbl[l]{BRoot}{0pt}{\lblsep}{\raisebox{-7pt}{\ensuremath{\scriptstyle 1}}}%
  \angled{AChild1}{AChild2}{2}%
  \angled{AChild2}{AChild3}{1}%
  \angled{BChild1}{BChild2}{2}%
}\ ,
\end{equation*}
$t_1$ verifies conditions~\refcond{InternalFertility}
and~\refcond{RootLabel} but does not verify condition~\refcond{IntermediateChildren},
$t_2$ verifies conditions~\refcond{InternalFertility} and~\refcond{IntermediateChildren} but
not~\refcond{RootLabel}, and $t_3$ verifies all three conditions.

The subindices in $\tii$ refer to the conditions imposed on the generator and on the \rb\ map of the free \rba, and they will be made clear in
Section~\ref{S:freeRB}. 

We define three subsets of $\tii$.  Let $\tiw$ be the subset of $\tii$
consisting of those trees whose internal node labels are less than or
equal to $1$. These elements can be seen as trees whose root label is $0$ or $1$ and the only other decorations are in the angles, since
the only possible label for the non-root internal nodes is $1$.  Let
$\twi$ be the subset of $\tii$ consisting of those trees whose angle
labels are $1$. These elements can be seen as trees whose only
decorations are on the internal nodes. Let $\tww = \tiw\cap\twi$.  The
set $\tww$ consists of two copies of (undecorated) rooted planar trees
satisfying conditions~\refcond{InternalFertility}
and~\refcond{IntermediateChildren}, where the label $0$ or $1$ at
the root of a tree indicates to which copy it belongs.
Table~\ref{fourSets} summarizes the decoration rules for each of the
four sets.

\begin{table}[!ht]
\centering
{\setlength\extrarowheight{3pt}
\begin{tabular}{|>{$}c<{$}||>{$}c<{$}|>{$}c<{$}|>{$}c<{$}|}
\hline
\text{Set} & \text{Root}  & \text{Angles} &
\vertical{\vborder{5pt}{\parbox{80pt}{\centering
  Non-root internal nodes}}}\\
\hhline{|=::===|}
\tii & \N & \Z^+ & \Z^+ \\
\tiw & \{0,1\}  & \Z^+ & \text{\{1\}}\\
\twi & \N & \text{\{1\}} & \Z^+  \\
\rule[-5pt]{0pt}{0pt}
\tww & \{0,1\}  & \text{\{1\}} & \text{\{1\}}\\
\hline
\end{tabular}}
\caption{Sets of decorated trees} \label{fourSets}
\end{table}

The following are examples of each kind of tree:
\begin{equation*} \label{examplesoftrees}
\begin{gathered}
  \vertical{\children{A}{,xbbh=\lblheight,xbbr=.4\treesepone,xbbd=\levelsep}{3}{1.2\treesepone}{\levelsep}%
    \rput[t](AChild2){%
      \children{B}{}{2}{.8\treesepone}{\levelsep}}%
    \rput[t](AChild3){%
      \children{C}{}{2}{.8\treesepone}{\levelsep}}%
    \lbl[bl]{ARoot}{0pt}{\lblsep}{2}%
    \lbl[r]{BRoot}{0pt}{\lblsep}{\raisebox{-7pt}{\ensuremath{\scriptstyle 3}}}%
    \lbl[l]{CRoot}{0pt}{\lblsep}{1}%
    \angled{AChild1}{AChild2}{4}%
    \angled{AChild2}{AChild3}{1}%
    \angled{BChild1}{BChild2}{2}%
    \angled{CChild1}{CChild2}{5}%
  }
  \in \tii,\quad
%%%
  \vertical{\children{A}{,xbbh=\lblheight,xbbr=.4\treesepone,xbbd=\levelsep}{3}{1.2\treesepone}{\levelsep}%
    \rput[t](AChild2){%
      \children{B}{}{2}{.8\treesepone}{\levelsep}}%
    \rput[t](AChild3){%
      \children{C}{}{2}{.8\treesepone}{\levelsep}}%
    \lbl[bl]{ARoot}{0pt}{\lblsep}{0}%
    \angled{AChild1}{AChild2}{4}%
    \angled{AChild2}{AChild3}{1}%
    \angled{BChild1}{BChild2}{2}%
    \angled{CChild1}{CChild2}{5}%
  }
  \in \tiw, \quad
%%%
  \vertical{\children{A}{,xbbh=\lblheight,xbbr=.4\treesepone,xbbd=\levelsep}{3}{1.2\treesepone}{\levelsep}%
    \rput[t](AChild2){%
      \children{B}{}{2}{.8\treesepone}{\levelsep}}%
    \rput[t](AChild3){%
      \children{C}{}{2}{.8\treesepone}{\levelsep}}%
    \lbl[bl]{ARoot}{0pt}{\lblsep}{2}%
    \lbl[r]{BRoot}{0pt}{\lblsep}{3}%
    \lbl[l]{CRoot}{0pt}{\lblsep}{1}%
  }
  \in \twi,\quad
%%%
  \vertical{\children{A}{,xbbh=\lblheight,xbbr=.4\treesepone,xbbd=\levelsep}{3}{1.2\treesepone}{\levelsep}%
    \rput[t](AChild2){%
      \children{B}{}{2}{.8\treesepone}{\levelsep}}%
    \rput[t](AChild3){%
      \children{C}{}{2}{.8\treesepone}{\levelsep}}%
    \lbl[bl]{ARoot}{0pt}{\lblsep}{1}%
  }
  \in \tww.
\end{gathered}
\end{equation*}

%When we later perform operations on the labels of decorated trees $t\in\txy$, the convention
%\begin{equation}\label{E:convention}
%1+1=1
%\end{equation}
%will be adopted,
%whenever we are adding angle labels and $i=2$ or whenever we are adding node labels and $j=2$.

%\smallskip

%The sets $\txy$ may also be seen as {\em quotients} of $\tii$, as indicated
%in the following diagram.
%\begin{equation} \label{D:four-spaces}
%  \begin{gathered}
%    {\psset{nodesep=1pt}%
%      \begin{psmatrix}[rowsep=8pt,colsep=10pt]
%        & \tii \\
%        \tiw & & \twi \\
%        & \tww \ncline{->>}{1,2}{2,3} \ncline{->>}{1,2}{2,1}
%        \ncline{->>}{2,1}{3,2} \ncline{->>}{2,3}{3,2}
%      \end{psmatrix}%
%    }%
%  \end{gathered}
%\end{equation}
%For example, the
%map $\twi\onto\tww$ sends a tree $t$ to a tree $t'$ with the same shape
%and with all positive node labels in $t$ changed into $1$ in $t'$. The other maps 
%are
%defined similarly, by replacing either angle labels or internal node labels $\geq 1$ by $1$, according to whether we change $i$ or $j$.
%The algebraic significance of these maps is explained in Remark~\ref{R:can-maps}.

% and the following diagram displays the inclusions among
%them.
%\begin{equation} \label{four-spaces}
%  {\psset{nodesep=2pt}
%    \setlength{\extrarowheight}{10pt}
%  \begin{array}{ccc}
%    & \rnode{A}{\tii} \\ 
%    \rnode{B}{\tiw} && \rnode{C}{\twi} \\
%    & \rnode{D}{\tww}
%  \end{array}
%  \ncline{A}{B}
%  \ncline{A}{C}
%  \ncline{B}{D}
%  \ncline{C}{D}
%}
%\end{equation}

\medskip

We consider two notions of degree for each kind of tree. The {\em node
  degree} of a decorated tree $t$ is the sum of the labels on the
internal nodes of $t$, and we denote it by $\degnode(t)$. Similarly,
the {\em angle degree} of $t$, denoted $\degangle(t)$, is the sum of
the labels in the angles of $t$.  Note that $\degangle(t)$ is always a
positive integer, while $\degnode$ may take the value $0$, namely, for
the trees $\Y{0}{i}$.

In particular, observe that for a tree $t$ in $\T_{2,j}$, the angle
degree coincides with the number of angles of $t$, which is one less than the number of leaves of $t$.  On the other hand, for
$t\in\T_{i,2}$, if the root of $t$ is labeled
by $1$ then the node degree coincides with the number of internal
nodes, while if it is labeled by $0$, the node degree is the number of
non-root internal nodes.

For $i,j\in\{2,\infty\}$, $n\ge 1$, and
$m\ge0$, let
\begin{equation*}
  \txy(n,m) = \bigl\{t\in\txy\mid \degangle(t)=n\text{\ and\ }
  \degnode(t)=m \bigr\}.
\end{equation*}
These sets will be linear bases for the homogeneous components of the
free \rba s, see Section~\ref{S:comb-dim}. In Table~\ref{T:listing} we
show the elements of $\tww(n,m)$ for $n=1,2,3$ and $m=0,1,2,3$.  We
set
\begin{equation*}
  \txy(*,m) = \bigsqcup_{n\ge 1}\txy(n,m),\quad
  \txy(n,*) = \bigsqcup_{m\ge 0}\txy(n,m),\quad
  \txy(k) = \bigsqcup_{\substack{n\ge 1,\;m\ge 0\\ n+m=k}} \txy(n,m).
\end{equation*}

\begin{table}[!ht]
  \centering
  \begin{tabular}{|>{$}c<{$}|>{$}c<{$}||>{$}c<{$}|}
    \hline
    n & m & \vborder{4pt}{\text{Elements of $\tww(n,m)$}} \\
    \hhline{|==::=|}
    1 & 0 & \vborder{4pt}{\ensuremath{\preYY{}{0}{}}} \\
    \hline
    1 & 1 & \vborder{4pt}{\ensuremath{\preYY{}{1}{}}} \\
    \hhline{|==::=|}
    2 & 0 & \vborder{4pt}{\text{\small empty}} \\
    \hline
    2 & 1 & \vertical{\includegraphics{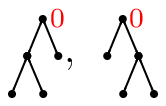}} \\
    \hline
    2 & 2 & \vertical{\includegraphics{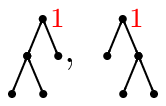}} \\
    \hhline{|==::=|}
    3 & 0 & \vborder{4pt}{\text{\small empty}} \\
    \hline
    3 & 1 & \vertical{\includegraphics{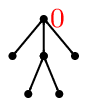}} \\
    \hline
    3 & 2 & \ \vertical{\includegraphics{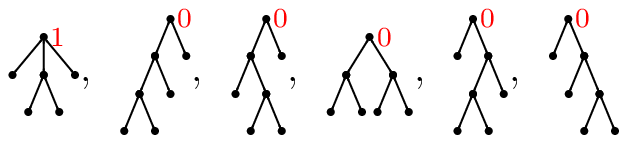}} \\
    \hline
    3 & 3 & \vertical{\includegraphics{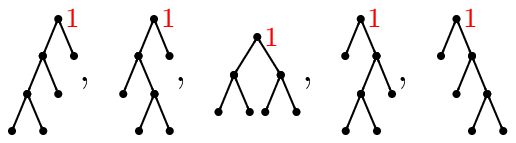}} \\
    \hline
  \end{tabular}
  \caption{$\tww(n,m)$ for $n=1,2,3$, and $m=0,1,2,3$}
  \label{T:listing}
\end{table}

We let $\txy^{+}$ (respectively, $\txy^{0}$) denote the subset of $\txy$ consisting of those trees whose root label is positive (respectively, $0$), and define $\txy^{a}(n,m)=\txy^{a}\cap\txy(n,m)$, for $a\in\{0,+\}$.

Let $ \widehat{\T}_{i,j} = \txy\cup\{\point\}$ be the set of
decorated trees with the (unlabeled) tree with a single node adjoined. 
We set $\deg(\point)=(0,0)$. Similarly, let  $\widehat{\T}^{a}_{i,j}= \txy^{a}\cup\{\point\}$ for $a\in\{0,+\}$.

\subsection{Grafting of decorated trees}\label{S:grafting}

We introduce a {\em grafting operation} on the set of decorated trees.   
Define a function
\begin{equation*}
 \Grxy :\bigcup_{n\ge 1} {(\widehat{\T}^+_{i,j})}^n\times{(\Z^+)}^{n-1} \longrightarrow \ \widehat{\T}_{i,j}
\end{equation*}
as follows. 
First, identify the set $\widehat{\T}^+_{i,j}\times {(\Z^+)}^0$ with $\widehat{\T}^+_{i,j}$ and set
\begin{equation*}
  \Grxy(t) = t. \label{E:def-G-1}
\end{equation*}
Then, for $n\geq 2$, set
\begin{equation}
  \Grxy(t_1,\ldots,t_n;\, i_1,\ldots, i_{n-1}) = N \biggl(
    \vertical{\childrenX{A}{,xbbh=5pt,xbbd=10pt,xbbr=4pt,xbbl=4pt}{4}{30pt}{5ex}%
      \lbl[bl]{ARoot}{0pt}{\lblsep}{0}%
      \lbl[t]{AChild1}{\lblsep}{0pt}{\textcolor{black}{%
          \ensuremath{\scriptstyle{\mathstrut t_1}}}}%
      \lbl[t]{AChild2}{\lblsep}{0pt}{\textcolor{black}{%
          \ensuremath{\scriptstyle{\mathstrut t_2}}}}%
      \lbl[t]{AChild3}{\lblsep}{0pt}{\textcolor{black}{%
          \ensuremath{\scriptstyle{\mathstrut t_{n-1}}}}}%
      \lbl[t]{AChild4}{\lblsep}{0pt}{\textcolor{black}{%
          \ensuremath{\scriptstyle{\mathstrut t_n}}}}%
      \angledG{AChild1}{AChild2}{\;\;i_1}%
      \angledG{AChild3}{AChild4}{i_{n-1}\;\;\:}%
      \angledX{AChild2}{AChild3}{{}\cdots{}}%
      \rput[t](ARoot|AChild2){\vborder{\lblsep}{\ensuremath{\scriptstyle{\mathstrut {}\cdots{}}}}}%
    }
    \biggr). \label{E:def-G-many}
\end{equation}
Here, the function $N$ normalizes the tree in such a way that the
result satisfies condition~\refcond{IntermediateChildren}; namely, if
$t_k=\point$, for $1<k<n$, then $t_k$ and the edge joining it to the
new root are removed from the tree, and the two adjacent angles (the
angle between $t_{k-1}$ and $t_k$ and the one between $t_k$ and
$t_{k+1}$) are merged into one angle which acquires the label
$i_{k-1}+i_k$.  Several additions may occur.

Another clarification is needed. When $i=2$, this addition is performed according to the convention 
\begin{equation}\label{E:convention-i}
1+1=1,
\end{equation}
so all angle labels remain equal to $1$. (Alternatively, if we view
trees in $\T_{2,j}$ as having no angle labels, then no additions are
necessary.)

In other words, for $n>1$, the operation $\Grxy$ grafts the trees
$t_k$ to a new root with label $0$, and uses the arguments $i_k$ as
the labels of the resulting new angles. Some of these are then added
if an intermediate leaf is formed. The result then satisfies
conditions~\refcond{InternalFertility},~\refcond{IntermediateChildren},
and~\refcond{RootLabel}, so it is a well-defined element of $\txy$.

For example,
\begin{equation*}
  \Grii\bigl(\betagenerator, \point, \point;\, 2, 1\bigr) = N\biggl(%
  \tree{\children{A}{,xbbh=\lblheight,xbbl=.4\treesepone,xbbd=\levelsep}{3}{\treesepone}{\levelsep}%
    \rput[t](AChild1){\children{B}{}{2}{.8\treesepone}{\levelsep}}%
    \lbl[bl]{ARoot}{0pt}{\lblsep}{0}%
    \lbl[r]{BRoot}{0pt}{\lblsep}{1}%
    \angled{AChild1}{AChild2}{2}%
    \angled{AChild2}{AChild3}{1}%
    \angled{BChild1}{BChild2}{1}%
  }
  \biggr) = 
  \tree{\children{A}{,xbbh=\lblheight,xbbl=.4\treesepone,xbbd=\levelsep}{2}{\treesepone}{\levelsep}%
    \rput[t](AChild1){\children{B}{}{2}{.8\treesepone}{\levelsep}}%
    \lbl[bl]{ARoot}{0pt}{\lblsep}{0}%
    \lbl[r]{BRoot}{0pt}{\lblsep}{1}%
    \angled{AChild1}{AChild2}{3}%
    \angled{BChild1}{BChild2}{1}%
  }.
\end{equation*}

We also define a {\em de-grafting operation} $\H_{i,j}:\txy \to \bigcup_{n\ge 1} {(\widehat{\T}^+_{i,j})}^n\times{(\Z^+)}^{n-1} $ by
\begin{equation} \label{E:def-H}
\H_{i,j}(t) = 
\begin{cases}
  (t_1,\ldots,t_n;\,  
   i_1,\ldots,i_{n-1}) & \text{if $t\in\txy^{0}$,} \\
  t & \text{if $t\in\txy^{+}$,}
\end{cases}
\end{equation}
where for $1\leq k\leq n$, $t_k$ is the subtree of $t$ rooted at the
$k$-th child of the root of $t$ (counting from left to right), and for
$1\leq k\leq n-1$, $i_k$ is the label of the angle between the $k$-th
and the $(k+1)$-th children.

For example,
\begin{equation*}
  \H_{\infty,\infty}\biggl( 
  \tree{\children{A}{,xbbh=\lblheight,xbbl=.4\treesepone,xbbd=\levelsep}{2}{\treesepone}{\levelsep}%
    \rput[t](AChild1){\children{B}{}{2}{.8\treesepone}{\levelsep}}%
    \lbl[bl]{ARoot}{0pt}{\lblsep}{0}%
    \lbl[r]{BRoot}{0pt}{\lblsep}{1}%
    \angled{AChild1}{AChild2}{3}%
    \angled{BChild1}{BChild2}{1}%
  }
  \biggr) = \bigl(\betagenerator, \point;\, 3\bigr),
  \quad\text{while}\quad
  \H_{\infty,\infty}\biggl(
  \tree{\children{A}{,xbbh=\lblheight,xbbl=.4\treesepone,xbbd=\levelsep}{2}{\treesepone}{\levelsep}%
    \rput[t](AChild1){\children{B}{}{2}{.8\treesepone}{\levelsep}}%
    \lbl[bl]{ARoot}{0pt}{\lblsep}{1}%
    \lbl[r]{BRoot}{0pt}{\lblsep}{1}%
    \angled{AChild1}{AChild2}{3}%
    \angled{BChild1}{BChild2}{1}%
  }
  \biggr) =
  \tree{\children{A}{,xbbh=\lblheight,xbbl=.4\treesepone,xbbd=\levelsep}{2}{\treesepone}{\levelsep}%
    \rput[t](AChild1){\children{B}{}{2}{.8\treesepone}{\levelsep}}%
    \lbl[bl]{ARoot}{0pt}{\lblsep}{1}%
    \lbl[r]{BRoot}{0pt}{\lblsep}{1}%
    \angled{AChild1}{AChild2}{3}%
    \angled{BChild1}{BChild2}{1}%
  }.
\end{equation*}

\subsection{Construction of the free \rba s on one generator}
\label{S:freeRB}

Let $\Rii^\lambda$ be the category
whose objects are triples $(A,x,\beta)$ where $(A,\beta)$ is a \rba\ 
and  $x\in A$ is an element.   A
morphism $f$ in $\Rii^\lambda$ from $(A,x,\beta)$ to $(B,y,\gamma)$ is a
 morphism of algebras that preserves the distinguished elements and
commutes with the \rb\ operators, that is,
\begin{equation*}
  f(x)=y,\qquad f\beta=\gamma f.
\end{equation*}

For $i,j\in\{2,\infty\}$, define $\Rxy^\lambda$ as the full
subcategory of $\Rii^\lambda$ whose objects $(A,x,\beta)$ satisfy that 
\[x^2=x \text{ \ if $i=2$, \ and \ } \beta^2=-\lambda\beta \text{ \ if $j=2$.}\]

By the {\em free \rba\ on one generator} we mean the initial object in the
category $\Rii^\lambda$. The initial object in $\Rwi^\lambda$ is the free \rba\ on one 
idempotent generator, the  initial object in $\Riw^\lambda$ is the free \rba\ on 
one generator and with a quasi-idempotent \rb\ map, and that in $\Rww^\lambda$ 
is the free  \rba\ on one idempotent generator and with a quasi-idempotent 
\rb\ map. 

The free \rba\  (the initial object in the category $\Rii^\lambda$)
was constructed by Ebrahimi-Fard and 
Guo~\cite{EG}. The  free \rba\ on one 
idempotent generator is also a special case of the constructions of~\cite{EG}.
 Below we provide a  simpler
description of these algebras, as well as of the related algebras mentioned
in the preceding paragraph, in a unified manner.

%\begin{equation*}
%\overline{\mathstrut\point}=\betaxy(\point)=\point,
%\end{equation*}
%and for $t\in\txy^+$ (respectively, $t\in\txy$),
%\begin{equation*}\label{E:beta}
% \overline{%
%   \psset{dotsize=0pt 1}%
%   \preYY{}{\mathstrut a}{}%
%   \rput[c](ARoot|IAChild1){\ensuremath{\scriptstyle{\textcolor{black}{t}}}}%
%   \ncline{AChild1}{AChild2}%
% }%
% =
% {%
%   \psset{dotsize=0pt 1}%
%   \preYY{,xbbr=12pt}{a-1}{}%
%   \rput[c](ARoot|IAChild1){\ensuremath{\scriptstyle{\textcolor{black}{t}}}}%
%   \ncline{AChild1}{AChild2}%
% }%
% \quad \text{and}\quad
% \betaxy\bigl({%
%   \psset{dotsize=0pt 1}%
%   \preYY{}{\mathstrut a}{}%
%   \rput[c](ARoot|IAChild1){\ensuremath{\scriptstyle{\textcolor{black}{t}}}}%
%   \ncline{AChild1}{AChild2}%
% }\bigr)%
% =
% {%
%   \psset{dotsize=0pt 1}%
%   \preYY{,xbbr=10pt}{a+1}{}%
%   \rput[c](ARoot|IAChild1){\ensuremath{\scriptstyle{\textcolor{black}{t}}}}%
%   \ncline{AChild1}{AChild2}%
% }.%
%\end{equation*}

\begin{definition}\label{D:product}
Fix $\lambda\in\field$. 
Let $\Fxy$ the vector space with basis $\txy$
and $\widehat{\F}_{i,j}$ the vector space with basis $\widehat{\T}_{i,j}$.
We extend the map  $\Grxy$ of Section~\ref{S:grafting} linearly to
these spaces. We define the map $\betaxy:\Fxy\to\Fxy$ as the linear
extension of
\begin{equation} \label{E:beta-map}
 \betaxy\bigl({%
   \psset{dotsize=0pt 1}%
   \preYY{}{\mathstrut a}{}%
   \rput[c](ARoot|IAChild1){\ensuremath{\scriptstyle{\textcolor{black}{t}}}}%
   \ncline{AChild1}{AChild2}%
 }\bigr)%
 =
 \begin{cases}
 {%
   \psset{dotsize=0pt 1}%
   \preYY{,xbbr=10pt}{a+1}{}%
   \rput[c](ARoot|IAChild1){\ensuremath{\scriptstyle{\textcolor{black}{t}}}}%
   \ncline{AChild1}{AChild2}%
}, & \text{when $j\not=2$;} \\[1ex]
 (-\lambda)^a   
 {%
   \psset{dotsize=0pt 1}%
   \preYY{,xbbr=0pt}{1}{}%
   \rput[c](ARoot|IAChild1){\ensuremath{\scriptstyle{\textcolor{black}{t}}}}%
   \ncline{AChild1}{AChild2}%
}, & \text{when $j=2$.}
 \end{cases}
\end{equation}
We also define $\betaxy(\point)=\point$ to extend the map to
$\betaxy:\widehat{\F}_{i,j}\to \widehat{\F}_{i,j}$.

We define a product $\astlambda$ on the space $\widehat{\F}_{i,j}$ and  a product $\assocprodlambda$ on the space $\Fxy$
by means of a mixed recursion. The recursion starts with
\begin{gather}\label{E:base-ast} 
  \point\astlambda u = u \astlambda \point = u \\
\intertext{for $u\in\widehat{\T}_{i,j}$, and follows with}
  \label{E:def-prod}
  t\assocprodlambda s  = 
 \Grxy \bigl( t_1,\ldots, t_{n-1}, \betaxy(\overline{t}_n
 \astlambda \overline{s}_1),  s_2, \ldots, s_m; \,
     i_1,\ldots, i_{n-1}, j_1,\ldots,j_{m-1} \bigr),\\
  \intertext{for $t$ and $s$ in $\txy$, and}
 \label{E:def-prodast}
   u\astlambda v = \betaxy(u)\assocprodlambda v + 
                  u \assocprodlambda \betaxy(v) + 
                  \lambda u\assocprodlambda v,
\end{gather}
for $u,v\in\widehat{\T}_{i,j}$. Here, we have set 
\begin{equation*}
\H(t) = (t_1,\ldots,t_n;\, i_1,\ldots,i_{n-1}) \quad\text{and}\quad
\H(s) = (s_1,\ldots,s_m;\, j_1,\ldots,j_{m-1})\,,
\end{equation*}
and $\overline{t}_n$ and $\overline{s}_1$ are
the result of the operation
\begin{equation*}
  \overline{{%
   \psset{dotsize=0pt 1}%
   \preYY{}{\mathstrut a}{}%
   \rput[c](ARoot|IAChild1){\ensuremath{\scriptstyle{\textcolor{black}{t}}}}%
   \ncline{AChild1}{AChild2}%
 }} =
  \begin{cases}
  {%
   \psset{dotsize=0pt 1}%
   \preYY{,xbbr=10pt}{a-1}{}%
   \rput[c](ARoot|IAChild1){\ensuremath{\scriptstyle{\textcolor{black}{t}}}}%
   \ncline{AChild1}{AChild2}%
}, & \text{if $a>0$;} \\[.5ex]
 \point, & \text{if $t=\point$}.
  \end{cases}
\end{equation*}
\end{definition}

%\begin{equation}
%  \begin{gathered}
%  \label{E:def-prod}
%  t\assocprodlambda s =\Grxy ( 
%  \begin{aligned}[t]
%    & t_1,\ldots, t_{n-1}, t_n \astlambda s_1,
%%    \betaxy(t_n)\assocprodlambda s_1 + t_n\assocprodlambda\betaxy(s_1)
%%    + \lambda t_n\assocprodlambda s_1, 
%    s_2, \ldots, s_m; \\
%    & i_1,\ldots, i_{n-1}, j_1,\ldots,j_{m-1} ),
%  \end{aligned}
%\end{gathered}
%\end{equation}

Note that $\overline{t}$ is undefined if the root label of $t$ is
$0$. In~\eqref{E:def-prod}, both $t_n$ and $s_1$ belong to
$\widehat{\T}^+_{i,j}$, so $\overline{t}_n$ and $\overline{s}_1$ are well
defined. In addition,
 $\overline{t}_n\astlambda \overline{s}_1$ involves the computation of products of the form
$t'\assocprodlambda s'$ satisfying $\degnode(t') \le \degnode(t)$ and
$\degnode(s')\le \degnode(s)$ with at least one of the inequalities
being strict. Thus~\eqref{E:def-prod} and~\eqref{E:def-prodast} invoke
each other recursively until either $t_n=\point$ or $s_1=\point$, at
which point the recursion stops with an application
of~\eqref{E:base-ast}.  In equation~\eqref{E:def-prod} we may
encounter a case when $n=1$ (or $m=1$). In such a case we understand
that the sequence $t_1,\ldots,t_{n-1}$ (or $s_2,\ldots,s_{m}$) is
empty, as usual.

By construction, the product $\astlambda$ is related to the product $\assocprodlambda$ and the operator $\betaxy$ by means of~\eqref{E:rb}. It will then follow, once we show that $(\Fxy,\assocprodlambda,\betaxy)$ is a \rba, that 
$(\Fxy,\astlambda)$ is an associative algebra, with $(\widehat{\F}_{i,j},\astlambda)$ being its unital augmentation (and with $\point$ being the unit element). Note, however, that the product
$\assocprodlambda$ is not defined on $\widehat{\F}_{i,j}$ and this space is not
a \rba.

\begin{example}\label{Ex:def-prod} We illustrate the definition of the product
$\assocprodlambda$ with a few small examples.
We have $\H\bigl(\Y{0}{i}\bigr) = (\point, \point;\, i)$. Therefore,
\begin{equation}\label{E:geni-genj}
 \Y{0}{i} \assocprodlambda \Y{0}{j} = 
    \Grii\bigl(\point,  \betaii(\point\astlambda \point), \point;\, i, j\bigr) 
  =\Grii(\point, \point, \point;\, i, j) = 
  \preYY{,treesep=20pt,levelsep=3.5ex}{0}{}%
  \ncline[linestyle=none]{ARoot}{AChild1}\ncput[npos=0.8]{\pnode{X}}%
  \ncline[linestyle=none]{ARoot}{AChild2}\ncput[npos=0.8]{\pnode{Y}}%
  \ncarc[linewidth=0.3pt,linecolor=blue,arcangle=-50,nodesep=1.5pt]{-}%
  {X}{Y}\naput[labelsep=1pt]%
  {\ensuremath{\scriptscriptstyle{\textcolor{blue}{i+j}}}}
\end{equation}
Also, since $\H\bigl(\Y{1}{i} \bigr)= \Y{1}{i}$, we have
\begin{equation}
\Y{1}{i} \assocprodlambda  \Y{0}{j} = \Grii\bigl( \betaii\bigl( \Y{0}{i} \astlambda
\point\bigr), \point;\, j\bigr) =\Grii\bigl( \Y{1}{i}, \point;\, j \bigr) =
  \tree{\children{A}{,xbbh=\lblheight,xbbl=.4\treesepone,xbbd=\levelsep}{2}{\treesepone}{\levelsep}%
    \rput[t](AChild1){\children{B}{}{2}{.8\treesepone}{\levelsep}}%
    \lbl[bl]{ARoot}{0pt}{\lblsep}{0}%
    \lbl[r]{BRoot}{0pt}{\lblsep}{1}%
    \angled{AChild1}{AChild2}{j}%
    \angled{BChild1}{BChild2}{i}%
  },
\end{equation}
using that $\betaii\bigl(\Y{0}{i}\bigr)=\Y{1}{i}$. With the same
considerations we obtain
\begin{align*}
\Y{0}{i}\astlambda  \Y{0}{j} = 
\Y{1}{i}\assocprodlambda  \Y{0}{j} +  \Y{0}{i} \assocprodlambda
\Y{1}{j} + \lambda \Y{0}{i}\assocprodlambda  \Y{0}{j} =
  \tree{\children{A}{,xbbh=\lblheight,xbbl=.4\treesepone,xbbd=\levelsep}{2}{\treesepone}{\levelsep}%
    \rput[t](AChild1){\children{B}{}{2}{.8\treesepone}{\levelsep}}%
    \lbl[bl]{ARoot}{0pt}{\lblsep}{0}%
    \lbl[r]{BRoot}{0pt}{\lblsep}{1}%
    \angled{AChild1}{AChild2}{j}%
    \angled{BChild1}{BChild2}{i}%
  }
  +
  \tree{\children{A}{,xbbh=\lblheight,xbbr=.4\treesepone,xbbd=\levelsep}{2}{\treesepone}{\levelsep}%
    \rput[t](AChild2){\children{B}{}{2}{.8\treesepone}{\levelsep}}%
    \lbl[bl]{ARoot}{0pt}{\lblsep}{0}%
    \lbl[l]{BRoot}{0pt}{\lblsep}{1}%
    \angled{AChild1}{AChild2}{i}%
    \angled{BChild1}{BChild2}{j}%
  }
  + \lambda
  \preYY{,treesep=20pt,levelsep=3.5ex}{0}{}%
  \ncline[linestyle=none]{ARoot}{AChild1}\ncput[npos=0.8]{\pnode{X}}%
  \ncline[linestyle=none]{ARoot}{AChild2}\ncput[npos=0.8]{\pnode{Y}}%
  \ncarc[linewidth=0.3pt,linecolor=blue,arcangle=-50,nodesep=1.5pt]{-}%
  {X}{Y}\naput[labelsep=1pt]%
  {\ensuremath{\scriptscriptstyle{\textcolor{blue}{i+j}}}}.
\end{align*}
Finally,
\begin{equation*}  
 \Y{1}{i} \assocprodlambda  \Y{1}{j} =\Grii\Bigl(\betaii\bigl(\Y{0}{i} \astlambda \Y{0}{i}\bigr)\Bigr)
=
  \tree{\children{A}{,xbbh=\lblheight,xbbl=.4\treesepone,xbbd=\levelsep}{2}{\treesepone}{\levelsep}%
    \rput[t](AChild1){\children{B}{}{2}{.8\treesepone}{\levelsep}}%
    \lbl[bl]{ARoot}{0pt}{\lblsep}{1}%
    \lbl[r]{BRoot}{0pt}{\lblsep}{1}%
    \angled{AChild1}{AChild2}{j}%
    \angled{BChild1}{BChild2}{i}%
  }
  +
  \tree{\children{A}{,xbbh=\lblheight,xbbr=.4\treesepone,xbbd=\levelsep}{2}{\treesepone}{\levelsep}%
    \rput[t](AChild2){\children{B}{}{2}{.8\treesepone}{\levelsep}}%
    \lbl[bl]{ARoot}{0pt}{\lblsep}{1}%
    \lbl[l]{BRoot}{0pt}{\lblsep}{1}%
    \angled{AChild1}{AChild2}{i}%
    \angled{BChild1}{BChild2}{j}%
  }
  + \lambda  
  \preYY{,treesep=20pt,levelsep=3.5ex}{1}{}%
  \ncline[linestyle=none]{ARoot}{AChild1}\ncput[npos=0.8]{\pnode{X}}%
  \ncline[linestyle=none]{ARoot}{AChild2}\ncput[npos=0.8]{\pnode{Y}}%
  \ncarc[linewidth=0.3pt,linecolor=blue,arcangle=-50,nodesep=1.5pt]{-}%
  {X}{Y}\naput[labelsep=1pt]%
  {\ensuremath{\scriptscriptstyle{\textcolor{blue}{i+j}}}}.
%  \preY{,treesep=18pt,levelsep=3.5ex}{1}{%
%    \raisebox{-10pt}{\ensuremath{\scriptscriptstyle{\textcolor{blue}{i+j}}}}}.
\end{equation*}
\end{example}

\medskip

Let $\Fxy^\lambda$ denote the space $\Fxy$ endowed with the product 
$\assocprodlambda$.

\begin{proposition} \label{mainthm}
  The initial object in the category $\Rxy^{\lambda}$ is $\bigl(\Fxy^\lambda,
  \generator, \beta_{i,j}\bigr)$.
\end{proposition}

\begin{proof} We first consider the case of $\Rii^\lambda$. This case is
  dealt with at length in~\cite{EG}, though in a
  different language.  We provide an independent proof
  %most of the required verifications,
  to illustrate the efficiency of the notation introduced in this paper.
  Our arguments extend to cover all categories $\Rxy^{\lambda}$, as  discussed at the end of the proof. 
  
During the course of the proof we omit the subindices from the symbols
$\Grii$, $\Hii$,
$\tii$, $\Fii$, $\widehat{\F}_{\infty,\infty}$, $\Rii$, and $\betaii$.
Thus, we abbreviate $\Gr = \Grii$, $\H = \Hii$, etc.
 We also fix de-grafting decompositions 
of trees $t$ and $s$~\eqref{E:def-H} as follows:
\begin{equation}\label{E:graf-dec}
\H(t) = (t_1,\ldots,t_n; \, i_1,\ldots, i_{n-1}), \quad 
\H(s)= (s_1,\ldots, s_m;\, j_1, \ldots, j_{m-1}).
\end{equation}

 We first check that $\beta$ is a \rb\ map. For any
 $t\in\T$, the root label of $\beta(t)$ is at least $1$, so $\beta(t)\in\T^+$ and
 by~\eqref{E:def-H} we have $\H\bigl(\beta(t)\bigr) =
 \beta(t)$. Using~\eqref{E:def-prod} and~\eqref{E:def-prodast} we obtain
 \begin{equation} \label{E:rb-map}
 \begin{aligned}
   \beta(t)\assocprodlambda \beta(s) &= \Gr \bigl( \beta(t\astlambda s )\bigr) =
   \beta( t\astlambda s)  \\
   &= \beta\bigl( 
       \beta(t)\assocprodlambda s + 
       t \assocprodlambda \beta(s) + 
       \lambda t\assocprodlambda s \bigr),
 \end{aligned}
\end{equation}
observing that $\overline{\beta(t)} = t$, by definition.
Hence, $\beta$ verifies condition~\eqref{E:rb} and it is a \rb~operator.
  
  Let $(A,x,\gamma)$ be an object of $\R^\lambda$.  
  We formally adjoin two elements $1_\ast$ and $1$ to $A$. We declare that $1_\ast$ is a unit element for the product $\ast$ of $A$~\eqref{E:rb}, and $1$ is a unit  element for the given product of $A$, and set
  \begin{equation*}
   \gamma(1_\ast)=1.
 \end{equation*}
In order to define a map
  $\varphi:\F\to A$, we first set
  \begin{equation*}
    \varphi(\point)=1_\ast,
  \end{equation*}
 and then, given $t\in\T$, define $\varphi(t)\in A$ recursively by
  \begin{equation} \label{E:def-phi}
    \varphi(t)= \gamma\bigl( \varphi(\overline{t}_1) \bigr) x^{i_1} 
                \gamma\bigl( \varphi(\overline{t}_2) \bigr) \cdots
                x^{i_{n-1}} \gamma\bigl( \varphi(\overline{t}_n) \bigr),
  \end{equation}
where $t_k$ and $i_k$ are as in~\eqref{E:graf-dec}.  When $t\in\T^+$, we have $\H(t)=t$, so this definition reads $\varphi(t)=\gamma\bigl( \varphi(\overline{t}) \bigr)$. Applying this to $t=\beta(s)$, where $s\in\T$ is an arbitrary tree, we obtain
  \begin{equation} \label{E:phi-commute}
    \varphi \bigl( \beta(s) \bigr) = \gamma\bigl(\varphi(s)\bigr),
  \end{equation}
since $\overline{t}=s$ in this case. Thus, $\varphi$ commutes with the \rb\ operators.
  
  Next, for $t=\generator$ we have
  $\H(t)=(\point,\point;\, 1)$ and hence
  \begin{equation*}
    \varphi\bigl(\generator\bigr) = \gamma\bigl(\varphi(\point) \bigr) x^1
                          \gamma\bigl(\varphi(\point) \bigr) =
    \gamma(1_\ast) x \gamma(1_\ast) = x,
  \end{equation*}
  proving that $\varphi$ preserves the distinguished elements of $\F$ and
  $A$.
  
  We now check that $\varphi$ is a morphism of algebras by induction
  on the bidegrees of $t$ and $s$. We will show that $\varphi$
  transforms the products $t\assocprodlambda s$ of $\F$ and
  $t\astlambda s$ of $\widehat{\F}$ into the products
  $\varphi(t)\varphi(s)$ of $A$ and
  $\varphi(t)\ast\varphi(s)$ of $\widehat{A}$, where $\widehat{A}$ is
  $A$ with the unit element $1_\ast$ adjoined. Since $\varphi$
  commutes with the \rb\ operators, knowing that
  $\varphi(t\assocprodlambda s)=\varphi(t)\varphi(s)$ holds up to a
  certain degree, implies that $\varphi(t\astlambda
  s)=\varphi(t)\ast\varphi(s)$ holds up to the same degree.  This will
  in turn be used to prove the former equality for the next degree.

The base case for the induction occurs when $t=\point$ or $s=\point$, for which it holds trivially that $\varphi(t\astlambda s)=\varphi(t)\ast\varphi(s)$.
  
 Before proceeding with the inductive step, we make a general observation.
  For a tuple of trees $u=(u_1,\ldots,u_n)\in(\widehat{\T}^+)^n$ and a tuple
  of positive integers $\alpha=(i_1,\ldots, i_{n-1})\in(\Z^+)^{n-1}$,
  we have
  \begin{equation} \label{E:phi-graft}
    \varphi\bigl(\Gr(u; \alpha)\bigr)
    =
    \gamma\bigl(\varphi(\overline{u}_1)\bigr) x^{i_1} \cdots
    x^{i_{n-1}} \gamma\bigl(\varphi(\overline{u}_n)\bigr).
  \end{equation}
  Indeed, if for all $k=2,\ldots,n-1$ we have $u_k\neq \point$, then $\H \bigl(
  \Gr(u; \alpha) \bigr) = (u;\alpha)$ and~\eqref{E:phi-graft} boils down to~\eqref{E:def-phi}. On the other hand, assume
  that $u_k=\point$ for exactly one $k=2,\ldots,n-1$. Then, 
  \begin{equation*}
    \H \bigl( \Gr(u; \alpha) \bigr) = 
    (u_1,\ldots,u_{k-1},u_{k+1},\ldots,u_n; \, i_1,\ldots,i_{k-1}+i_k,\ldots, i_{n-1})
  \end{equation*}
  and therefore
  \begin{equation*}
    \varphi\bigl(\Gr(u;\alpha) \bigr) =
    \gamma\bigl(\varphi(\overline{u}_1)\bigr) x^{i_1} \cdots x^{i_{k-1}+i_k}
    \cdots x^{i_{n-1}}\gamma\bigl(\varphi(\overline{u}_n)\bigr),
  \end{equation*}
  but as $x^{i_{k-1}+i_k} =
  x^{i_{k-1}}\gamma\bigl(\varphi(\point)\bigr) x^{i_k}$, we see
  that~\eqref{E:phi-graft} holds for such $u$. A similar argument applies if more than one $u_k=\point$, and we conclude
  that~\eqref{E:phi-graft} holds for any   tuple $u\in(\widehat{\T}^+)^n$.

  Applying equation~\eqref{E:phi-graft} to the tuple entering in  definition~\eqref{E:def-prod}, we obtain
  \begin{equation} \label{E:phi-prod}
    \varphi(t\assocprodlambda s) = 
      \gamma\bigl(\varphi(\overline{t}_1)\bigr) x^{i_1} \cdots x^{i_{n-1}}
      \gamma\bigl(\varphi(\overline{t}_n\astlambda \overline{s}_1)\bigr) 
      x^{j_1} \cdots x^{j_{m-1}} \gamma\bigl(\varphi(\overline{s}_m)\bigr).
  \end{equation}
%  Consider the factor $\gamma\bigl(\varphi(t_n\astlambda
%  s_1)\bigr)$.  
%   If $t_n=\point$, then $t_n\astlambda s_1=s_1$ and the
%  right hand side is $\varphi(t)\varphi(s)$, by~\eqref{E:def-phi}. 
%  The same applies if $s_1=\point$. 
  The bidegrees of $\overline{t}_n$ and $\overline{s}_1$ are smaller than those of $t$ and $s$, so
we can assume inductively that $\varphi(\overline{t}_n\astlambda \overline{s}_1) =
  \varphi(\overline{t}_n)\astlambda \varphi(\overline{s}_1)$. Therefore,
  \begin{equation*}
    \gamma\bigl(\varphi(\overline{t}_n\astlambda \overline{s}_1) \bigr) = \gamma\bigl(
    \varphi(\overline{t}_n)\astlambda \varphi(\overline{s}_1) \bigr) =
    \varphi(\overline{t}_n)\varphi(\overline{s}_1),
  \end{equation*}
  the latter equality in view of~\eqref{E:rb}. Substituting 
  in~\eqref{E:phi-prod} and recalling~\eqref{E:def-phi} we obtain $\varphi(t\assocprodlambda s) =  \varphi(t)\varphi(s)$, which completes the induction.

\smallskip

It remains to prove that the product $\assocprodlambda$ is associative. 
Let $t$,
$s$, and $u\in\T$. In addition to the de-grafting decompositions of
$t$ and $s$%~\eqref{E:graf-dec}
, let $\H(u) = (u_1,\ldots, u_\ell;\,
k_1,\ldots,k_{\ell -1})$ be that of $u$. Two main cases arise, according to $m$, the number of trees in the  decomposition of $s$. If
$m=1$, we have 
\begin{equation*}
\H(t\assocprodlambda s) =\H\Gr \bigl(t_1,\ldots, t_{n-1},\beta(\overline{t}_n\astlambda
\overline{s});\, i_1,\ldots,i_{n-1}\bigr)=\bigl(t_1,\ldots, t_{n-1},\beta(\overline{t}_n\astlambda
\overline{s});\, i_1,\ldots,i_{n-1}\bigr)
\end{equation*}
and
\begin{equation*}
\H(s\assocprodlambda u) =\H\Gr \bigl(\beta(\overline{s}\astlambda \overline{u}_1),u_2,\ldots,u_\ell;\, k_1,\ldots,k_{\ell-1}\bigr)=\bigl(\beta(\overline{s}\astlambda \overline{u}_1),u_2,\ldots,u_\ell;\, k_1,\ldots,k_{\ell-1}\bigr).
\end{equation*}
 since all intermediate trees in these tuples are different from $\point$. 
 % Do we need this remark here??
% Observe that the
%$\astlambda$ product only yields $\point$ when both factors are
%$\point$. 
Therefore,
\begin{equation*}\label{E:assoc-1}
  (t\assocprodlambda s)\assocprodlambda u = \Gr\bigl(t_1,\ldots,
  \beta\bigl((\overline{t}_n\astlambda \overline{s})\astlambda \overline{u}_1\bigr),\ldots, u_\ell;\,
  i_1,\ldots,i_{n-1},k_1, \ldots, k_{\ell-1}\bigr).
\end{equation*}
and
\begin{equation*}
  t\assocprodlambda (s\assocprodlambda u) = \Gr\bigl(t_1,\ldots,
  \beta\bigl(\overline{t}_n\astlambda (\overline{s}\astlambda \overline{u}_1)\bigr),\ldots, u_\ell;\,
  i_1,\ldots,i_{n-1},k_1, \ldots, k_{\ell-1}\bigr).
\end{equation*}
We can assume inductively on the node degree that $\assocprodlambda$
is associative on smaller degree trees, which implies that $\astlambda$
is also associative on those trees. Therefore,
$(\overline{t}_n\astlambda \overline{s} )\astlambda
\overline{u}_1=\overline{t}_n\astlambda (\overline{s}\astlambda
\overline{u}_1)$ and
$(t\assocprodlambda s)\assocprodlambda u= t\assocprodlambda(s\assocprodlambda u)$. The base case occurs when
some of the trees $t_n$, $s$ or $u_1$ equals $\point$, and then associativity holds trivially.

If $m>1$, then in the tuples
\begin{equation*}
\bigl(t_1,\ldots,\beta(\overline{t}_n\astlambda \overline{s}_1),\ldots,s_m\bigr) \quad\text{and}\quad
\bigl(s_1,\ldots,\beta(\overline{s}_m\astlambda \overline{u}_1),\ldots,u_\ell\bigr)
\end{equation*}
the only intermediate trees that may equal $\point$ are $\beta(\overline{t}_n\astlambda
\overline{s}_1)$ and $\beta(\overline{s}_m\astlambda \overline{u}_1)$. If neither of them equals $\point$
then
\begin{equation*}
\H(t\assocprodlambda s)=\bigl(t_1,\ldots,\beta(\overline{t}_n\astlambda \overline{s}_1),\ldots,s_m;\,i_1,\ldots,i_{n-1},j_1,\ldots,j_{m-1}\bigr)
\end{equation*}
and
\begin{equation*}
 \H(s\assocprodlambda u) =\bigl(s_1\ldots,\beta(\overline{s}_m\astlambda \overline{u}_1),\ldots,u_\ell;\,
  k_1,\ldots,k_{\ell-1}\bigr),
\end{equation*}
from which it follows that both $(t\assocprodlambda s)\assocprodlambda u$ and $t\assocprodlambda(s\assocprodlambda u)$ equal
\begin{equation*}
  \Gr\bigl(t_1,\ldots,\beta(\overline{t}_n\astlambda \overline{s}_1),\ldots,\beta(\overline{s}_m\astlambda \overline{u}_1),\ldots,u_\ell;\,
i_1,\ldots,i_{n-1},j_1,\ldots,j_{m-1},k_1,\ldots,k_{\ell-1}\bigr)\,.
\end{equation*}

It remains to deal with the cases when $\overline{t}_n\astlambda
\overline{s}_1=\point$ or $\overline{s}_m\astlambda \overline{u}_1=\point$. We consider the case
when $\overline{t}_n\astlambda \overline{s}_1=\point$ and $\overline{s}_m\astlambda
\overline{u}_1\not=\point$;  the others are similar.  We have,
\begin{equation*}
  \H(t\assocprodlambda s) = (t_1,\ldots, t_{n-1}, s_2, \ldots, s_m; \,i_1,
  \ldots, i_{n-1}+j_1, \ldots, j_{m-1})
\end{equation*}
and
\begin{equation*}
  \H(s\assocprodlambda u) = \bigl(s_1,\ldots, \beta(\overline{s}_m\astlambda \overline{u}_1),\ldots u_\ell;\,
  j_1,\ldots, j_{m-1}, u_1, \ldots, u_{\ell-1}\bigr).
\end{equation*}
Hence, both $ (t\assocprodlambda s)\assocprodlambda u$ and 
$t\assocprodlambda(s\assocprodlambda u)$ equal
\begin{equation*} 
   \Gr \bigl( t_1,\ldots,   t_{n-1}, s_2, \ldots, \beta(\overline{s}_m\astlambda \overline{u}_1), \ldots, u_\ell; \, i_1, \ldots, i_{n-1}+j_1,
    \ldots, j_{m-1}, k_1, \ldots , k_{\ell-1}\bigr).
\end{equation*}
This completes the proof of the proposition for the
case of $\Rii^\lambda$.

\smallskip

Most of the preceding proof goes through for the general case of the category $\Rxy^\lambda$. To finish the proof, we comment on the few exceptional situations that arise when $i=2$ or $j=2$. 
 
When $i=2$. First note that the element $\generator$ is indeed idempotent, in view of~\eqref{E:convention-i} and~\eqref{E:geni-genj}. Now,
in the proof of equation~\eqref{E:phi-graft} we encounter
$i_{k-1}+i_k = 1+1 = 1$. However, since $(A,x,\gamma)$ is an object of $\R_{2,j}^\lambda$, we have $x^2=x$, and therefore
$x$ can still be split as
$x=x\gamma\bigl(\varphi(\point) \bigr) x$. Thus
equation~\eqref{E:phi-graft} and  all the conclusions
about the map $\varphi$ are valid.

When $j=2$. First note that $\beta_{i,2}$ is quasi-idempotent: for any  $t\in\T_{i,2}$, the root label of $\beta_{i,2}(t)$ is $1$, so by~\eqref{E:beta-map},
$\beta_{i,2}\bigl(\beta_{i,2}(t)\bigr)=-\lambda\beta_{i,2}(t)$.
In relation to the case $j=\infty$,  only the proofs of~\eqref{E:rb-map}  and~\eqref{E:phi-commute} require additional argument. The reason is that, in view of~\eqref{E:beta-map},
 $\beta_{i,2}(t)$ depends on the root label of $t$, which
may be $0$ or $1$, and in the latter case we get $\beta_{i,2}(t) =
-\lambda t$. In the former case, $\overline{\beta_{i,2}(t)} = t$.

The proof of~\eqref{E:rb-map} still holds when $t$ and $s$ belong to $\T^0_{i,2}$. Suppose that $t\in\T^+_{i,2}$ and $s\in \T^0_{i,2}$ 
(the other cases are similar). In this case,
\begin{equation*}
\beta_{i,2}(t)\assocprodlambda\beta_{i,2}(s)=-\lambda t\assocprodlambda
\beta_{i,2}(s) =-\lambda\Gr_{i,2}\bigl(\beta_{i,2}(\overline{t}\astlambda s)\bigr)=-\lambda\beta_{i,2}(\overline{t}\astlambda s).
\end{equation*}
On the other hand, $\beta_{i,2}(t)\assocprodlambda s=-\lambda t\assocprodlambda s$, so
\begin{equation*}
\beta_{i,2}(t\astlambda s) =
\beta_{i,2}\bigl(t\assocprodlambda\beta_{i,2}(s)\bigr) =
\beta_{i,2}\Bigl(\Gr_{i,2}\bigl(\beta_{i,2}(\overline{t}\astlambda s)\bigr)\Bigr) =
\beta_{i,2}\bigl(\beta_{i,2}(\overline{t}\astlambda s)\bigr).
\end{equation*}
Since $\beta_{i,2}$ is quasi-idempotent, we conclude $\beta_{i,2}(t)\assocprodlambda\beta_{i,2}(s)=\beta_{i,2}(t\astlambda s)$, as needed.

The proof of~\eqref{E:phi-commute}  still holds if $s\in\T^0_{i,2}$.  
If $s\in\T^+_{i,2}$ then we have $
\beta_{i,2}(\overline{s})=s $ and $\beta_{i,2}(s)=-\lambda s$. Hence,
  \begin{equation*}
  \varphi \beta_{i,2}(s) =-\lambda\varphi(s) =-\lambda \varphi
  \beta_{i,2}(\overline{s}) = -\lambda\gamma 
  \varphi(\overline{s}) 
    = \gamma  \gamma
  \varphi(\overline{s}) =
  \gamma  \varphi
  \beta_{i,2}(\overline{s})  = \gamma \varphi(s).
\end{equation*}
%  \begin{multline*}
%  \varphi\bigl( \beta_{i,2}(s) \bigr)=-\lambda\varphi(s) =-\lambda \varphi\bigl(
%  \beta_{i,2}(\overline{s}) \bigr) = -\lambda\gamma \bigl(
%  \varphi(\overline{s}) \bigr) \\
%  = \gamma^2 \bigl(
%  \varphi(\overline{s}) \bigr) = \gamma \Bigl( \gamma\bigl(
%  \varphi(\overline{s}) \bigr) \Bigr) =
%  \gamma \Bigl( \varphi\bigl(
%  \beta_{i,2}(\overline{s}) \bigr) \Bigr) = \gamma \bigl( \varphi(s) \bigr).
%\end{multline*}
We used that  $\gamma$ is a quasi-idempotent \rb\ operator, which holds since in this case $(A,x,\gamma)$ is an object of $\R_{i,2}^\lambda$,
and that $\varphi \beta_{i,2}(\overline{s}) = \gamma 
  \varphi(\overline{s}) $, which holds since $\overline{s}\in\T^0_{i,2}$.
  This completes the proof of the proposition.
\end{proof}

\begin{remark}\label{R:adjoint} Consider the forgetful functor from the category
of \rba s to the category of algebras.
The adjoint functor was constructed by Ebrahimi-Fard and Guo~\cite{EG}.
Applying this functor to the one-dimensional algebra $\field\{x\}$ ($x^2=x$) yields the algebra $\Fwi^\lambda$, while applying it to the algebra
$x\field[x]$ yields $\Fii^\lambda$. Our notation is also useful for describing this functor: simply consider decorated trees in which the angle labels are elements of a given algebra $A$. The notion of grafting naturally extends to this context
(using the product of $A$ when a merging of angles occurs in~\eqref{E:def-G-many}), and the constructions of this section carry through. The result is the value of the adjoint functor on the algebra~$A$.
\end{remark}

We derive a useful recursive expression for the canonical morphism from the free \rba\ to another \rba.

\begin{corollary}\label{C:can-exp}
Let $(A,x,\beta)$ be an object  in $\Rxy^\lambda$, and $\varphi:\Fxy^\lambda\to A$  the unique morphism of \rba s such that $\varphi(\generator)=x$.  Given a tree $t\in\txy$, let $a\in\N$ be its root label, $t_1,\ldots,t_n\in\widehat{\T}_{i,j}^+$ be the subtrees of $t$ rooted at the children of the root of $t$, and $i_1,\ldots,i_{n-1}$  the labels of the angles between these children, as in~\eqref{E:def-H}. Then,
\begin{equation}\label{E:can-exp}
\varphi(t)= \gamma^a\bigl( \varphi(t_1)  x^{i_1} 
                 \varphi(t_2) x^{i_2}  \cdots
                x^{i_{n-1}}  \varphi(t_n)\bigr).
\end{equation}
In particular, $t$ decomposes as
\begin{equation}\label{E:can-dec}
t= \betaxy^a\Bigl( t_1\assocprodlambda \bigl( \generator\bigr)^{i_1} \assocprodlambda
    t_2\assocprodlambda\bigl(\generator\bigr)^{i_2} \assocprodlambda \cdots     \assocprodlambda  \bigl(\generator\bigr)^{i_{n-1}} \assocprodlambda t_n\Bigr),
\end{equation}
where we understand that if $t_1$ or $t_n$ are equal to $\point$ then they are omitted.
\end{corollary}
\begin{proof}
  As shown in the proof of Proposition~\ref{mainthm}, the map
  $\varphi$ is defined by~\eqref{E:def-phi}. If $a\geq 1$, then
  $\Hxy(t)=t$, so
  $\varphi(t)=\gamma\bigl(\varphi(\overline{t})\bigr)$. The root label
  of $\overline{t}$ is $a-1$. Proceeding by induction we see that
  $\varphi(t)=\gamma^a\bigl(\varphi(\hat{t})\bigr)$, where $\hat{t}$
  is the same tree as $t$ but with root label $0$. Now, since
  $\Hxy(\hat{t})=(t_1,\ldots,t_n;\,
  i_1,\ldots,i_{n-1})$, and
  $\gamma\bigl(\varphi(\overline{t}_k)\bigr)=\varphi\bigl(\betaxy(\overline{t}_k)\bigr)=\varphi(t_k)$,
  an application of~\eqref{E:def-phi} gives
\begin{equation*}
  \varphi(\hat{t})=\gamma\bigl( \varphi(\overline{t}_1) \bigr) x^{i_1} 
  \gamma\bigl( \varphi(\overline{t}_2) \bigr) \cdots
  x^{i_{n-1}} \gamma\bigl( \varphi(\overline{t}_n) \bigr)=
  \varphi(t_1)  x^{i_1} 
  \varphi(t_2) x^{i_2}  \cdots
  x^{i_{n-1}}  \varphi(t_n),
\end{equation*}
and ~\eqref{E:can-exp} follows.            

Applying this result to $A=\Fxy^\lambda$, $x=\generator$, and $\gamma=\betaxy$ we obtain~\eqref{E:can-exp}, since in this case $\varphi$ is the identity.
\end{proof}

  The inclusions among the various categories $\Rxy^\lambda$ determine
  morphisms in the opposite direction among the corresponding initial objects,
  as indicated below
\begin{equation*}
    \vertical{\ensuremath{\psset{nodesep=1pt}%
      \begin{psmatrix}[rowsep=8pt,colsep=10pt]
        & \Rii^\lambda \\
        \Riw^\lambda & & \Rwi^\lambda \\
        & \Rww^\lambda
        \ncline{-}{1,2}{2,3} \ncline{-}{1,2}{2,1}
        \ncline{-}{2,1}{3,2} \ncline{-}{2,3}{3,2}
      \end{psmatrix}%
    }}%
    \qquad
    \vertical{\ensuremath{\psset{nodesep=1pt}%
      \begin{psmatrix}[rowsep=8pt,colsep=10pt,labelsep=2pt]
        & \Fii^\lambda \\
        \Fiw^\lambda & & \Fwi^\lambda \\
        & \Fww^\lambda
        \ncline{->>}{1,2}{2,3}\naput{\scriptstyle{\varphi_{\cdot,\infty}}} 
        \ncline{->>}{1,2}{2,1}\nbput{\scriptstyle{\varphi_{\infty,\cdot}}}
        \ncline{->>}{2,1}{3,2}\nbput{\scriptstyle{\varphi_{\cdot,2}}} 
        \ncline{->>}{2,3}{3,2}\naput{\scriptstyle{\varphi_{2,\cdot}}}
      \end{psmatrix}%
    }}%
\end{equation*}
These maps are the unique morphisms of \rba s that preserve
the distinguished elements $\generator{}$. We describe these maps next.

\begin{proposition}\label{P:can-maps}
The maps $\varphi_{\cdot,j}$ are the linearizations of the maps
$\T_{\infty,j}\onto\T_{2,j}$ that erase all angle labels. 
The maps $\varphi_{i,\cdot}$ are given by
\begin{equation*}\label{E:quotient-map}
t\mapsto (-\lambda)^d t'
\end{equation*}
where the tree $t'$ is obtained from $t$ by changing
 all positive node labels into $1$, and the exponent $d$ is equal to
 $\degnode(t)$ minus the number of nodes of $t$ with positive labels.
\end{proposition}
\begin{proof} This follows from~\eqref{E:beta-map} and~\eqref{E:can-exp}.
\end{proof}

We conclude this section by discussing a canonical filtration on the free \rba s.
For each $a\geq 0$, let $\Fxy^{\lambda,a}$ be the subspace of $\Fxy^\lambda$
spanned by those trees $t\in\txy$ with root node label less than or equal to $a$.
Thus,
\[\Fxy^{\lambda,0}\subseteq \Fxy^{\lambda,1}\subseteq\Fxy^{\lambda,2}\subseteq\cdots\]
is an increasing sequence of subspaces of $\Fxy^{\lambda}$. Notice that
$\F_{i,2}^{\lambda,a}=\F_{i,2}^{\lambda,1}$ for any $a\geq 1$.

Recall (Definition~\ref{D:product}) that $\widehat{\F}_{i,j}^{\lambda}$ denotes
the unital augmentation of the algebra $\Fxy^{\lambda}$. 
Define $\widehat{\F}_{i,j}^{\lambda,a}$ as
the span of $\Fxy^{\lambda,a}$ and the single node $\point$. 
In particular, 
$\point\in \widehat{\F}_{i,j}^{\lambda,0}$. 

\begin{proposition}\label{P:filtration} For any $a,b\geq 0$, $\betaxy(\Fxy^{\lambda,a})\subseteq \Fxy^{\lambda,a+1}$, and
\[\widehat{\F}_{i,j}^{\lambda,a}\astlambda \widehat{\F}_{i,j}^{\lambda,b}\subseteq \widehat{\F}_{i,j}^{\lambda,a+b},
\text{ \ and \ }
\Fxy^{\lambda,a}\assocprodlambda\Fxy^{\lambda,b}\subseteq
\begin{cases}
\Fxy^{\lambda,a+b-1} & \text{ if $a>0$ and $b>0$,}\\
 \Fxy^{\lambda,0}        & \text{ if $a=0$ or $b=0$.}
\end{cases}\]
In particular, $\widehat{\F}_{i,j}^{\lambda,0}$ is a unital subalgebra of 
$(\widehat{\F}_{i,j}^{\lambda},\astlambda)$ and $ \Fxy^{\lambda,0}$ is an ideal of
$(\Fxy^{\lambda},\assocprodlambda)$.
\end{proposition}
\begin{proof} This may be proved by induction, using~\eqref{E:base-ast},~\eqref{E:def-prod}, and~\eqref{E:def-prodast}.
\end{proof}

\section{Combinatorics of free Baxter algebras}\label{S:combinatorics}

\subsection{Trees and paths}
\label{S:comb-trees}

We establish several bijections between the sets of trees defined in
Section~\ref{S:decorated} and other combinatorial objects of a more familiar nature. This is used in Section~\ref{S:comb-dim}
to compute the dimensions of the homogeneous
components of the free \rba s. The bijections are in the same spirit as those in~\cite[Proposition 6.2.1]{St99}.

Let us set the notation for the sets of combinatorial objects. Let
$\PT$ be the set of  rooted planar trees whose internal nodes have at least
two children. For $n\ge 1$ and $m\ge 0$, let $\PT(n,m)$ be the subset
of $\PT$ consisting of trees with $n+1$ leaves and $m$ internal nodes.
Also let $\PT(n)$ be the set of planar rooted trees with $n+1$ leaves,
so that
\begin{equation*}
  \PT(n) = \bigsqcup_{m\ge 0}\PT(n,m).
\end{equation*}
Observe that $\PT(n,0)=\emptyset$ for any $n\ge 1$. The
cardinality of $\PT(n)$ is the small Schr\"oder number~\cite[Exercise
6.39]{St99}. For $m\geq 1$, the cardinality of $\PT(n,m)$ is
\begin{equation}\label{E:plan-trees}
\frac{1}{n+1}\binom{n+m}{m}\binom{n-1}{m-1}.
\end{equation}
This is also the number of $(n-m)$-dimensional faces of the $(n-1)$-dimensional associahedron~\cite[Exercise 6.33]{St99}.

Let
$\BinT(n)=\PT(n,n)$. This is the subset of $\PT(n)$ consisting of binary trees.   Its cardinality  is the Catalan
number~\cite[Exercise 6.19]{St99}
\begin{equation*}
  C(n) = \frac{1}{n+1}\binom{2n}{n}.
\end{equation*}
This is also the number of vertices of the $(n-1)$-dimensional associahedron.

Next we define various sets of lattice paths.
\begin{descr}
\item[$\CP(n)$] set of  Catalan paths of length $2n$; that is,
  lattice paths from $(0,0)$ to $(n,n)$ with steps $H=(1,0)$ and
  $V=(0,1)$, never rising above the diagonal. The number of these
  paths is the Catalan number $C(n)$.
\item[$\SP(n)$] set of Schr\"oder paths of length $2n$; that is, lattice paths
  from $(0,0)$ to $(n,n)$ with steps $H=(1,0)$, $V=(0,1)$, and
  $D=(1,1)$, never rising above the diagonal.
The number of these paths is
  the large Schr\"oder number (twice the small Schr\"oder number)~\cite[Exercise 6.39]{St99}.
\item[$\SP(n,m)$] set of Schr\"oder paths of length $2n$ with $n-m$
    diagonal steps.  The number of these paths is given in Proposition~\ref{P:fiw}.
\item[$\RP(n,m)$] set of paths in $\SP(n,m)$ such that each diagonal
  step  is followed by a horizontal step, except if it is the last
  step. The number of these paths is given in Proposition~\ref{P:fww}.
\item[$\MP(n)$] set of  Motzkin paths of length $n$; that is,
  lattice paths from $(0,0)$ to $(n,0)$ with steps $U=(1,1)$,
  $H=(1,0)$, and $D=(1,-1)$, never crossing below the $x$-axis. The
  number of these paths is the Motzkin number~\cite[Exercise
  6.38]{St99}.
\item[$\RMP(n)$] set of paths in $\MP(n)$ such that each horizontal
  step is followed by an up step, except if it is the last step.
\item[$\hMP(n)$] set of $h$-colored Motzkin paths from $(0,0)$ to
  $(n-1,0)$; that is, Motzkin paths whose horizontal steps are colored
  with one of two colors. The number of these paths is
  $C(n)$~\cite[Exercise (yyy)]{catadd}.
\item[$\huMP(n)$] set of $(h,u)$-colored Motzkin paths from $(0,0)$
  to $(n-1,0)$; that is, Motzkin paths whose horizontal and up steps
  are colored with one of two colors. The number of these paths enters
  in Proposition~\ref{P:fww}.
 % This is sequence A071356.
\end{descr}
A few examples follow. The letters \textsf{R} and
\textsf{B} stand for the colors of the steps.
\begin{equation*}
  \setlength\extrarowheight{2ex}
  \begin{matrix}
  \vc{\begin{pspicture}(0,-0.1)(4,4)
    \psset{subgriddiv=0,gridlabels=0pt,dimen=middle,gridwidth=.3pt}
    \psline[linestyle=dotted](0,0)(4,4)
    \psgrid
    \psset{linewidth=1.5pt}
    \psline{-c}(0,0)(1,0)
    \psline{c-c}(1,0)(2,1)
    \psline{c-c}(2,1)(3,1)
    \psline{c-c}(3,1)(3,2)
    \psline{c-c}(3,2)(4,2)
    \psline{c-c}(4,2)(4,3)
    \psline{c-}(4,3)(4,4)
    \rput[r](-.3,0){$\scriptstyle (0,0)$}
    \rput[l](4.3,4){$\scriptstyle (4,4)$}
  \end{pspicture}}
  &\qquad\qquad &
  \vc{\begin{pspicture}(4,2)
    \psset{linewidth=0.5pt}
    \psline(0,0)(4,0)
    \psset{linewidth=1.5pt}
    \psline{c-c}(0,0)(1,0)
    \psline{c-c}(1,0)(2,1)
    \psline{c-c}(2,1)(3,1)
    \psline{c-c}(3,1)(4,0)
    \rput[r](-.3,0){$\scriptstyle 0$}
    \rput[l](4.3,0){$\scriptstyle 4$}
  \end{pspicture}}
  &\qquad\qquad &
  \vc{\begin{pspicture}(4,2)
    \psset{linewidth=0.5pt}
    \psline(0,0)(4,0)
    \psset{linewidth=1.5pt}
    \psline{c-c}(0,0)(1,0)
    \psline{c-c}(1,0)(2,1)
    \psline{c-c}(2,1)(3,1)
    \psline{c-c}(3,1)(4,0)
    \rput[r](-.3,0){$\scriptstyle 0$}
    \rput[l](4.3,0){$\scriptstyle 4$}
    \rput[c](0.5,0.4){\tiny\textsf{\textcolor{red}{R}}}
    \rput[c](1.3,0.9){\tiny\textsf{\textcolor{blue}{B}}}
    \rput[c](2.5,1.4){\tiny\textsf{\textcolor{blue}{B}}}
  \end{pspicture}} \\
  \RP(4,3) & & \MP(4) & & \huMP(5)
\end{matrix}
\end{equation*}

\medskip

The set of Schr\"oder paths $\SP(n,m)$ and its subset $\RP(n,m)$ can be
decomposed into two disjoint subsets:
\begin{equation*}
  \SP(n,m) = \SPp(n,m)\sqcup \SPz(n,m), \quad
  \RP(n,m) = \RPp(n,m)\sqcup \RPz(n,m),
\end{equation*}
where $\SPp(n,m)$ (respectively $\RPp(n,m)$) consists of those paths
in $\SP(n,m)$ (respectively $\RP(n,m)$) which do not have diagonal
steps lying on the diagonal, and $\SPz(n,m)$ (respectively
$\RPz(n,m)$) is its complement in $\SP(n,m)$ (respectively
$\RP(n,m)$).

\begin{proposition} \label{P:bijections}
  Let $n\ge 1$ and $m\ge 0$.
  \begin{enumerate}[(i)]
  \item \label{I:phi+} The sets $\tiw^+(n,m)$, $\PT(n,m)$, and
    $\SPp(n,m)$ are in bijection:
    \begin{equation*}
      \varphi^+: \tiw^+(n,m) \overset{f^+}{\longrightarrow} \PT(n,m) 
      \overset{g^+}{\longrightarrow} \SPp(n,m).
    \end{equation*}
  \item \label{I:phi0} The sets $\tiw^0(n,m)$, $\PT(n,m+1)$, and
    $\SPz(n,m)$ are in bijection:
    \begin{equation*}
      \varphi^0: \tiw^0(n,m) \overset{f^0}{\longrightarrow} \PT(n,m+1)
      \overset{g^0}{\longrightarrow} \SPz(n,m).
    \end{equation*}
    Moreover, there is a  bijection 
    $T:\SPp(n,m+1)\to \SPz(n,m)$ making the following diagram commutative, where $\betaiw$ is the map that changes the root label from $0$ to~$1$ (Section~\ref{S:grafting}),
    \begin{equation} \label{E:phi0-commute}
      \begin{gathered}
      {\psset{nodesep=2pt}%
        \begin{psmatrix}[rowsep=20pt,colsep=30pt]
          \tiw^0(n,m)  & \tiw^+(n,m+1) \\ 
          \SPz(n,m) & \SPp(n,m+1)
          \ncline{->}{1,1}{1,2}\naput{\scriptstyle{\betaiw}}
          \ncline{->}{1,1}{2,1}\nbput{\scriptstyle{\varphi^0}}
          \ncline{->}{1,2}{2,2}\naput{\scriptstyle{\varphi^+}}
          \ncline{->}{2,2}{2,1}\nbput{\scriptstyle{T}}
        \end{psmatrix}%
      }%
    \end{gathered}
  \end{equation}
  \item The bijections $ \varphi^+$ and $ \varphi^0$ 
  %in~$(\ref{I:phi+})$ and~$(\ref{I:phi0})$
    restrict to bijections
    \begin{equation*}
      \psi^+:\tww^+(n,m) \to \RP^+(n,m),\quad 
      \psi^0:\tww^0(n,m) \to \RP^0(n,m).
    \end{equation*}
\end{enumerate}
\end{proposition}

\begin{proof}[Proof of part (i)]
  Given a tree $t\in\tiw^+(n,m)$, define $f^+(t)$ as the planar tree
  resulting from substituting the decorations $j$ in each angle for
  $j-1$ intermediate leaves in the corresponding node. The tree
  $f^+(t)$ will have $n$ angles, hence $n+1$ leaves, and $m$ internal
  nodes, since $\degnode(t)$, for $t\in\twi^+$, coincides with the
  number of internal nodes. This process is clearly bijective. For
  example,
  \begin{equation*}
    \tree{\children{A}{,xbbd=\levelsep,xbbh=\lblheight}{3}{\treesepone}{\levelsep}%
      \rput[t](AChild2){%
        \children{B}{}{2}{\treesepone}{\levelsep}}%
      \lbl[bl]{ARoot}{0pt}{\lblsep}{1}%
      \angled{AChild1}{AChild2}{1}%
      \angled{AChild2}{AChild3}{2}%
      \angled{BChild1}{BChild2}{3}%
    }
    \quad\overset{f^+}{\longmapsto}\quad
    \tree{\children{A}{,xbbd=\levelsep,xbbl=.15\levelsep}{4}{.6\treesepone}{\levelsep}%
      \rput[t](AChild2){%
        \children{B}{}{4}{.5\treesepone}{\levelsep}}}
  \end{equation*}

  To define the function $g^+$, consider a tree $t\in\PT(n,m)$. We
  generate a Schr\"oder path $p\in\SPp(n,m)$ using
  Algorithm~\ref{treetopath}, see
  Appendix~\ref{S:appendix}. Informally, Algorithm~\ref{treetopath}
  traverses the tree depth-first and generates an $H$ step when it finds the
  leftmost child of a node, a $D$ step when it finds an intermediate
  child, and a $V$ step when it finds the rightmost child. For
  example,
  \begin{equation*}
    \tree{\children{A}{,xbbh=\lblheight,xbbr=.4\treesepone,xbbd=\levelsep}{3}{1.2\treesepone}{\levelsep}%
      \rput[t](AChild2){%
        \children{B}{}{2}{.8\treesepone}{\levelsep}}%
      \rput[t](AChild3){%
        \children{C}{}{2}{.8\treesepone}{\levelsep}}%
      \lbl[bl]{ARoot}{0pt}{\lblsep}{1}%
    }
    \quad\overset{g^+}{\longmapsto}\quad
    \vertical{\begin{pspicture}(0,-0.1)(4,4)
        \psset{subgriddiv=0,gridlabels=0pt,dimen=middle,gridwidth=.3pt}
        \psline[linestyle=dotted](0,0)(4,4)
        \psgrid
        \psset{linewidth=1.5pt}
        \psline{-c}(0,0)(1,0)
        \psline{c-c}(1,0)(2,1)
        \psline{c-c}(2,1)(3,1)
        \psline{c-c}(3,1)(3,2)
        \psline{c-c}(3,2)(3,3)
        \psline{c-c}(3,3)(4,3)
        \psline{c-}(4,3)(4,4)
      \end{pspicture}}    
  \end{equation*}
  
  The proof that Algorithm~\ref{treetopath} stops is straightforward,
  since the number of nodes of the trees involved in the recursive
  invocations of {\sffamily TreeToPath} is strictly less than that
  of~$t$.
  
  Note that for each internal node visited by the algorithm, an $H$
  step is issued when descending to its leftmost child. Similarly, for
  each angle a $V$ step or a $D$ step is issued when descending to
 an intermediate child or to the rightmost child.  Therefore, the path contains
  $m$ horizontal steps, and $n$ steps which are either vertical or
  diagonal. It can easily be proved by induction on the number of
  nodes that the algorithm generates an underdiagonal path with the
  same number of horizontal steps as of vertical steps. Therefore, the
  path must contain $m$ vertical steps and $n-m$ diagonal steps, and
  go from $(0,0)$ to $(n,n)$.
  
  We claim that a path generated by Algorithm~\ref{treetopath} cannot
  have a diagonal step lying on the diagonal. Suppose this were the
  case. The diagonal step cannot be the last step of the path, since
  Algorithm~\ref{treetopath} ends issuing a vertical step.  After
  issuing such a diagonal step, the algorithm processes the tree
  $t_i$. The result is a portion of the path that returns to the
  diagonal.  Eventually a vertical step is issued on exiting the
  innermost \textbf{if}, which would make the path cross  the
  diagonal. This proves that $g^+(t)\in\SPp(n,m)$.
  
  For the reverse process we use Algorithm~\ref{pathtotree} in
  Appendix~\ref{S:appendix}. {}From a path in $\SPp(n,m)$ we generate a
  tree, starting from a single node, by creating the children
  according to the steps of the path, read from $(0,0)$ to $(n,n)$. If
  the step is $H$, then a new child is created, the node is marked as
  available for creating more children, and the algorithm descends to
  the newly created child. If the step is $D$, then a new child is
  created in the first available node, searching upward from the
  current position. The same happens if the step is $V$, but in this
  case the node where the child is created is marked as no longer
  available.

  It is easy to see, inductively on $n$, that a path in $\SPp(n,m)$
  yields a planar tree. Note that the internal nodes are created by
  horizontal steps, hence there are $m$ such nodes. And the $D$ and
  $V$ steps produce angles, thus there are $(n-m)+m=n$
  angles. Therefore, the output of Algorithm~\ref{pathtotree} is a
  tree in $\PT(n,m)$. Clearly, the two algorithms are inverse of each
  other.
\end{proof}

\begin{proof}[Proof of part (ii)]
  We define the map $T:\SPp(n,m+1)\to\SPz(n,m)$ shown in
  Diagram~\eqref{E:phi0-commute} as follows. Let $p=s_1s_2\cdots s_k$
  be a path in $\SPp(n,m+1)$. Since a path in $\SPp(n,m+1)$ does not
  have $D$ steps on the main diagonal, $s_1$ and $s_k$ must be $H$ and
  $V$ steps, respectively. Let $p'$ be the path obtained from $p$ by dropping the frist and last steps, and shifting the rest by $(-1,0)$. In other words,
  $p'=s_2\cdots s_{k-1}$, with origin at $(0,0)$ and end at $(n-1,n-1)$.
  Define
  \begin{equation*}\label{E:T}
    T(p) = 
    \begin{cases}
      s_2\cdots s_{k-1}D & \text{if $p'$ is underdiagonal,}
      \\
      s_2\cdots s_{i-1}D H s_{i+1} \cdots s_{j-1}V
      s_{j+1} \cdots s_{k-1} &
      \text{if $p'$ crosses the diagonal,}
    \end{cases}
  \end{equation*}
  where, in the second case, $s_i$ is the first vertical step of $p'$ above
  the diagonal  and $s_j$ is the last step of $p'$ to the left of, or on,
  the diagonal. In both cases, the path $T(p)$ goes from $(0,0)$ to
  $(n,n)$,  is underdiagonal and has one more diagonal step than
  $p'$. Since the original path $p$ has $n-(m+1)$ diagonal steps,
  after applying $T$ we are left with a path with $n-(m+1)+1 = n-m$
  diagonal steps. It is also clear that 
  the path $T(p)$ always has a $D$ step on the main diagonal. Therefore,
  $T(p)\in\SPz(n,m)$.
  
  To define the inverse of $T$ observe that if $p'$ crosses the diagonal
  then $T(p)$ does not end in a $D$ step. If this were the case, $s_{k-1}=D$
  would be a step of $p'$ on the diagonal, which contradicts the choice of $j$.
   Therefore, the inverse of
  $T$ can be defined by
  \begin{equation*} \label{E:Tinverse}
    T^{-1}(p) = 
    \begin{cases}
      Hs_2\cdots s_{k-1}V & \text{if $s_k=D$,} \\
      Hs_1\cdots s_{i-1}V s_{i+2}\cdots s_{j-1} H s_{j+1}\cdots s_k V
      &
      \text{if $s_k\not= D$,}
    \end{cases}
  \end{equation*}
  where $s_i$ is the first diagonal step on the diagonal and $s_j$
  is the first vertical step after $s_i$ which touches the
  diagonal. 
  This proves that $T$ is bijective.
  
  Observe that the map $\betaiw:\tiw^0(n,m)\to\tiw^+(n,m+1)$ is a
  bijection. The bijection $\varphi^0$ is constructed as
  \begin{equation*}
    \tiw^0(n,m) \xrightarrow{\betaiw} \tiw^+(n,m+1) 
    \xrightarrow{f^+} \PT^+(n,m+1) \xrightarrow{g^+}
    \SPp(n,m+1) \xrightarrow{T} \SPz(n,m),
  \end{equation*}
  which fills Diagram~\eqref{E:phi0-commute}.
\end{proof}

\begin{proof}[Proof of (iii)]
  The function $f^+$ restricted to $\tww^+(n,m)$ just erases the root
  labels of the trees. Then, given a tree $t\in\tww^+$, we only need
  to verify that every $D$ step in $\varphi^+(t)$ is followed by an $H$
  step. But after Algorithm~\ref{treetopath} issues a $D$ step, the
  tree $t_i$ is processed and since it cannot be a leaf by
  condition~\refcond{IntermediateChildren} the next issued step must
  be $H$. Conversely, when processing a $D$ step followed by an $H$ step, 
  Algorithm~\ref{pathtotree} creates a child of the
  intermediate node we are visiting, and hence no intermediate node is
  a leaf. This proves that $\varphi^+$ restricts to a bijection
  $\psi^+:\tww^+(n,m)\to\RPp(n,m)$.
  
  Also, the function $T$ preserves the condition that $D$ steps are
  followed by $H$ steps, then Diagram~\eqref{E:phi0-commute}
  shows that $\varphi^0$ restricts to a bijection
  $\psi^0:\tww^0(n,m)\to\RPz(n,m)$.
\end{proof}

\begin{remark}\label{R:BSS} The bijections of Proposition~\ref{P:bijections} give a description for the number of
$(n-m)$-dimensional faces of the $(n-1)$-dimensional associahedron~\eqref{E:plan-trees} in terms of two classes of Schr\"oder paths (the sets $\SP^+(n,m)$ and $\SP^0(n,m-1)$). A description in
terms of a different class of Schr\"oder paths is given in~\cite[Proposition 2.7]{BSS}.
\end{remark}

\begin{corollary} \label{C:bijections}
  \begin{enumerate}[(i)]
  \item For $n\ge1$ and $m\ge 0$, there are bijections
    \begin{gather*}
      \tiw(n,m) \leftrightarrow \PT(n,m)\sqcup\PT(n,m+1)
      \leftrightarrow \SP(n,m), \\
      \tww(n,m) \leftrightarrow \RP(n,m).
    \end{gather*}
  \item For $n\ge 1$, there are bijections
    \begin{gather*}
      \tiw^+(n,\ast) \leftrightarrow \tiw^0(n,\ast) \leftrightarrow
      \PT(n) \leftrightarrow \SPp(n) \leftrightarrow \SPz(n), \\
      \tiw(n,\ast) \leftrightarrow \PT(n)\times\{0,+\} \leftrightarrow
      \SP(n). 
    \end{gather*}
    In particular,  $\#\SPp(n)=\#\SPz(n)=$ small
    Schr\"oder number.
  \item For $n\ge 1$, there are bijections
    \begin{gather}
        \tww^+(n,\ast) \leftrightarrow \tww^0(n,\ast) \leftrightarrow
        \huMP(n) \label{E:huMP(n)} \\
      \tww(n,\ast) \leftrightarrow \huMP(n)\times\{0,+\}. \label{E:2huMP(n)}
  %\huMP(n)\sqcup\huMP(n)
    \end{gather}
  \end{enumerate}
\end{corollary}

\begin{proof}
  Parts $(i)$ and $(ii)$ are immediate from
  Proposition~\ref{P:bijections}.  For part $(iii)$, we construct the
  bijection $f:\tww^+(n,\ast)\to \huMP(n)$ as follows. Given
  $t\in\tww^+$, consider the path $\varphi^+(t) = s_1s_2\cdots
  s_k\in\RPp(n,m)$ for some $m$. We know that $s_1=H$ and $s_k=V$.
  Consider the path $p=s_2\cdots s_{k-1}$ and start reading it from left
  to right. Using Table~\ref{T:substitutions}, the first time that one
  of the patterns listed in the left column of the table is found,
  write the corresponding value of the right column, and continue with
  the rest of the path.
\begin{table}[!ht]
\centering
\begin{threeparttable}[b]
\begin{tabular}{c|c|c|c}
% Line 1
\hhline{~|--|~}
\hspace{100pt} & 
\rule[-7pt]{0pt}{20pt}Pattern & Substitution \\
\hhline{~|==|~}
&
\vertical{\begin{pspicture}(0,-0.5)(1,1.5)
\psline{*-*}(0,0)(1,0)
\psline{*-*}(1,0)(1,1)
\end{pspicture}}
&
\vertical{\begin{pspicture}(0,0)(1,1)
\psline{*-*}(0,0)(1,0)
\rput[c](0.5,0.5){\tiny\textsf{\textcolor{red}{R}}}
\end{pspicture}}
& \hspace{100pt}
\\
\hhline{~|--|~} &
% Line 7
\vertical{\begin{pspicture}(0,-0.5)(3,1.5)
\psline{*-*}(0,0)(1,0)
\psline{*-*}(1,0)(2,1)
\psline{*-*}(2,1)(3,1)
\end{pspicture}}
&
\vertical{\begin{pspicture}(0,-0.5)(2,2.2)
\psline{*-*}(0,0)(1,1)
\psline{*-*}(1,1)(2,1)
\rput[c](0.3,0.9){\tiny\textsf{\textcolor{blue}{B}}}
\rput[c](1.5,1.5){\tiny\textsf{\textcolor{red}{R}}}
\end{pspicture}}
\\
\hhline{~|--|~} &
% Line 3
\vertical{\begin{pspicture}(0,-0.5)(2,0.5)
\psline{*-*}(0,0)(1,0)
\psline{*-*}(1,0)(2,0)
\end{pspicture}}
&
\vertical{\begin{pspicture}(0,-0.5)(1,1.5)
\psline{*-*}(0,0)(1,1)
\rput[c](0.3,0.9){\tiny\textsf{\textcolor{red}{R}}}
\end{pspicture}}
\\
\hhline{~|--|~} &
% Line 5
\vertical{\begin{pspicture}(0,-0.5)(2,2.5)
\psline{*-*}(0,0)(1,1)
\psline{*-*}(1,1)(2,1)
\psline{*-*}(2,1)(2,2)
\end{pspicture}}
&
\vertical{\begin{pspicture}(0,0)(2,1.4)
\psline{*-*}(0,0)(1,1)
\psline{*-*}(1,1)(2,0)
\rput[c](0.3,0.9){\tiny\textsf{\textcolor{blue}{B}}}
\end{pspicture}}
\\
\hhline{~|--|~} &
% Line 8
\vertical{\begin{pspicture}(0,-0.5)(4,2.5)
\psline{*-*}(0,0)(1,1)
\psline{*-*}(1,1)(2,1)
\psline{*-*}(2,1)(3,2)
\psline{*-*}(3,2)(4,2)
\end{pspicture}}
&
\vertical{\begin{pspicture}(0,-0.5)(3,2.5)
\psline{*-*}(0,0)(1,1)
\psline{*-*}(1,1)(2,2)
\psline{*-*}(2,2)(3,1)
\rput[c](0.2,0.8){\tiny\textsf{\textcolor{blue}{B}}}
\rput[c](1.2,1.8){\tiny\textsf{\textcolor{blue}{B}}}
\end{pspicture}}
\\
\hhline{~|--|~} &
% Line 6
\vertical{\begin{pspicture}(0,-0.5)(3,1.5)
\psline{*-*}(0,0)(1,1)
\psline{*-*}(1,1)(2,1)
\psline{*-*}(2,1)(3,1)
\end{pspicture}}
&
\vertical{\begin{pspicture}(0,-0.5)(2,2.2)
\psline{*-*}(0,0)(1,1)
\psline{*-*}(1,1)(2,1)
\rput[c](0.3,0.9){\tiny\textsf{\textcolor{blue}{B}}}
\rput[c](1.5,1.5){\tiny\textsf{\textcolor{blue}{B}}}
\end{pspicture}}
\\
\hhline{~|--|~} &
% Line 4
\vertical{\begin{pspicture}(0,-0.5)(0,2.5)
\psline{*-*}(0,0)(0,1)
\psline{*-*}(0,1)(0,2)
\end{pspicture}}
&
\vertical{\begin{pspicture}(0,0)(1,1.4)
\psline{*-*}(1,0)(0,1)
\end{pspicture}}
\\
\hhline{~|--|~} &
% Line 9
\vertical{\begin{pspicture}(0,-0.5)(2,2.5)
\psline{*-*}(0,0)(0,1)
\psline{*-*}(0,1)(1,2)
\psline{*-*}(1,2)(2,2)
\end{pspicture}}
&
\vertical{\begin{pspicture}(0,0)(1,1)
\psline{*-*}(1,0)(0,1)
\end{pspicture}}
\tnote{(*)}
\\
\hhline{~|--|~} &
% Line 2
\vertical{\begin{pspicture}(0,-0.5)(1,1.5)
\psline{*-*}(0,0)(0,1)
\psline{*-*}(0,1)(1,1)
\end{pspicture}}
&
\vertical{\begin{pspicture}(0,0)(1,1)
\psline{*-*}(0,0)(1,0)
\rput[c](0.5,0.5){\tiny\textsf{\textcolor{blue}{B}}}
\end{pspicture}}
\\
\hhline{~|--|~} 
\end{tabular}
\begin{tablenotes}%[para]
\item[(*)] also search backwards the rightmost point where
  the Motzkin path up-crossed the current level
%    \begin{pspicture}[.3](-.2,0)(1.1,1)
%      \psline{*-*}(0,0)(1,1)
%      \rput[c](0.2,0.8){\tiny\textsf{R}}
%    \end{pspicture}
    and insert before it the step 
    \begin{pspicture}[.3](-.2,0)(1.1,1)
      \psline{*-*}(0,0)(1,1)
      \rput[c](0.2,0.8){\tiny\textsf{\textcolor{blue}{B}}}
    \end{pspicture}.
\end{tablenotes}
\caption{Conversion to colored Motzkin paths.}
\label{T:substitutions}
\end{threeparttable}
\end{table}
Let $p'$ be the resulting path.
Consider the increment in the distance to the
diagonal, from the start to the end point, for each pattern of $p$.
When this increment is $0$, so is the increment of distance to the
line $y=0$ in the path $p'$. Note that line 8 in 
Table~\ref{T:substitutions} is one of these cases. The condition that
$\varphi^+(t)$ does not have diagonal steps lying on the diagonal
guarantees that the end point of the down step in $p'$ is above the
line $y=0$, and thus, the search indicated in line 8 is not empty.
 When the increment in $p$ is $\pm \sqrt{2}$
(diagonally), the increment in $p'$ is $\pm 1$ (vertically),
in each case with the same sign. Observe that for each pattern of $p$ of
length $2k$, the corresponding portion of the path $p'$ has length
$k$. Hence, the path $p'$ goes from $(0,0)$ to $(n-1,0)$, as $p$ has
length $2n-2$. Moreover, after removing the first and last steps of
$\varphi^+(t)$,  $p$ rises above the diagonal by at
most $\sqrt{2}/2$ (diagonally). Since this difference is not enough for the
path $p'$ to cross below the horizontal line $y=0$, by the previous
argument, we conclude that $p'$ is a Motzkin path in $\huMP(n)$.

For the reverse process it is enough to use
Table~\ref{T:substitutions} from right to left, taking into account
that patterns that are not in the right column (for example a
\textsf{B}-up step followed by a \textsf{R}-up step) come from an
application of line 8. After adding an $H$ step at the beginning of the
resulting Schr\"oder path, and a $V$ step at the end, we guarantee
that the result is in $\RPp(n,m)$ for some $m$. {}From there use the
inverse of $\psi^+$ to get a tree in $\tww^+(n,\ast)$.
\end{proof}

\begin{remark}
  The bijection~\eqref{E:huMP(n)} in Corollary~\ref{C:bijections} is
  an extension of the bijection $\CP(n) \leftrightarrow \hMP(n)$
  proposed by Stanley as solution to Exercise (yyy) in~\cite{catadd}.
  More precisely, consider the function $i:\hMP(n)\hookrightarrow
  \huMP(n)$ which sends a Motzkin path in which only the horizontal steps
are  colored, to the same Motzkin path with all up steps colored
   red (\textsf{R}). It is easy to see by looking at
  Table~\ref{T:substitutions} that the corresponding Schr\"oder path
  for $i(p)$, $p\in\hMP(n)$, under the bijection $(iii)$ in
  Corollary~\ref{C:bijections}, is actually a Catalan path. If we
  embed $\BinT(n)$ in $\tww^+(n,\ast)$ as trees with root label
  $1$, then the bijection $\psi^+$ from $(iii)$ in
  Proposition~\ref{P:bijections} also restricts and yields the
  commutative diagram
  \begin{equation*}
    {\psset{nodesep=3pt}%
      \begin{psmatrix}[rowsep=20pt,colsep=30pt]
        \tww^+(n,\ast) & \RPp(n,\ast) & \huMP(n)  \\ 
        \rule{0pt}{12pt} \BinT(n) & \rule{0pt}{12pt} \CP(n) & 
        \rule{0pt}{12pt} \hMP(n)
        \ncline{<->}{1,1}{1,2}
        \ncline{<->}{1,2}{1,3}
        \ncline{<->}{2,1}{2,2}
        \ncline{<->}{2,2}{2,3}
        \ncline{->}{2,1}{1,1}\ncput[npos=0]{\pnode{A}}%
        \ncline{->}{2,2}{1,2}\ncput[npos=0]{\pnode{B}}%
        \ncline{->}{2,3}{1,3}\ncput[npos=0]{\pnode{C}}%
        \rput(A){\pnode(-.4,0){AA}}\ncarc[arcangle=90,ncurv=1,nodesep=0pt]{A}{AA}
        \rput(B){\pnode(-.4,0){BB}}\ncarc[arcangle=90,ncurv=1,nodesep=0pt]{B}{BB}
        \rput(C){\pnode(-.4,0){CC}}\ncarc[arcangle=90,ncurv=1,nodesep=0pt]{C}{CC}
      \end{psmatrix}%
    }%
  \end{equation*}
\end{remark}

To close this section we study the subspaces $\tiw(k)$ and
$\tww(k)$. In this cases, too, it is possible to construct bijections
with familiar combinatorial objects, such as Motzkin paths. Recall that
\begin{equation*}
  \txy(k) = \bigsqcup_{\substack{n\ge 1, m\ge 0\\ n+m=k}} \txy(n,m).
\end{equation*}

\begin{proposition}\label{P:fiw-k}
For $k\ge 1$, there are bijections
\begin{equation*}
%  \bigsqcup_{\substack{n\ge 1, m\ge 0\\ n+m=k}} \tiw(n,m) 
  \tiw(k)
  \leftrightarrow \MP(k),
  \quad
%  \bigsqcup_{\substack{n\ge 1, m\ge 0\\ n+m=k}} \tww(n,m) 
  \tww(k)
  \leftrightarrow \RMP(k).
\end{equation*}
\end{proposition}

\begin{proof}
  {}From a tree $t\in\tiw(n,m)$, with $n+m=k$, we get the path
  $p\in\SP(n,m)$ using part $(i)$ of Corollary~\ref{C:bijections}.
  Reflect the path over the diagonal and rotate it clockwise until the
  diagonal becomes horizontal. If we consider each step to be of
  length $1$, we obtain a Motzkin path whose length is
  $2m+(n-m)=n+m=k$, since $p$ has $m$ horizontal steps, $m$ vertical
  steps, and $n-m$ diagonal steps.  For the converse, given a Motzkin
  path $p$ of length $k$, let $m$ be the number of up steps, which
  must coincide with the number of down steps. Let $n=k-m$. Then, the
  number of horizontal steps is $k-2m=n-m$. After rotating
  counter-clockwise and reflecting the path along the diagonal we obtain
  a Schr\"oder path $q$. This path has $m$ horizontal steps, $m$
  vertical steps, and $n-m$ diagonal steps. Hence, we conclude that
  $q\in\SP(n,m)$ with $n+m=k$. Apply again the bijection of
  Corollary~\ref{C:bijections} to obtain a tree in $\tiw(n,m)$.
  
  The second bijection is just the restriction of the previous one.
  Indeed, the bijection $(i)$ in~\ref{C:bijections} restricts to
  $\tww(n,m) \to \RP(n,m)$ and the condition about diagonal steps
  followed by horizontal steps translates, after the geometric
  transformations, to a condition about horizontal steps followed by
  up steps.
\end{proof}

\subsection{Bigrading and dimensions of the homogeneous components}
\label{S:comb-dim}

Recall the functions $\degangle:\txy\to\Z^+$ and $\degnode:\txy\to\N$ defined in Section~\ref{S:decorated}. Consider the bigrading on the vector space $\Fxy$ defined by
 \begin{equation}\label{E:deg-R}
\deg(t)=\bigl(\degangle(t),\degnode(t)\bigr).
\end{equation}
Let $\Fxy(n,m)$ denote  the homogeneous component
of bidegree $(n,m)$, so that the set $\txy(n,m)$ is
a basis for $\Fxy(n,m)$ and
\begin{equation*}
  \Fxy = \bigoplus_{\substack{n\ge 1\\ m\ge 0}} \Fxy(n,m).
\end{equation*}
Similarly, let 
\begin{equation*}
\Fxy(n,*) = \bigoplus_{m\ge 0} \Fxy(n,m), \quad
\Fxy(*,m)=  \bigoplus_{n\ge 1} \Fxy(n,m),
\end{equation*}
and 
\begin{equation*}
\Fxy(k)= \bigoplus_{\substack{n\ge 1,\, m\ge 0 \\n+m=k}}\Fxy(n,m).
\end{equation*}
Thus, $\Fxy(n,*)$, $\Fxy(*,m)$, and $\Fxy(k)$ are the subspaces of
$\Fxy$ spanned by $\txy(n,*)$, $\txy(*,m)$, and
$\txy(k)$, respectively. 

%Formula~\eqref{E:deg-R} can be understood as follows:
% $\deg(t)=(n,m)$ if when writing $t\in\F_{i,j}$ as a \emph{reduced expression} in the symbols  $x=\generator$ and $\beta_{i,j}$ with respect to the product $\assocprodlambda$, the symbol $x$ occurs exactly $n$ times and the symbol $\betaxy$ occurs \emph{at least} $m$ times 
% (due to the non-homogeneity of~\eqref{E:rb},
% there may be different ways to express $t$ in the symbols  $x$ and $\beta_{\infty,j}$, and
 
 The following assertions may be proven by induction from Definition~\ref{D:product}.
 If $i=\infty$, the product $\assocprodlambda$
preserves the angle degree. If $i=2$ (so $x$ is idempotent), the angle degree of a  $\assocprodlambda$ product  may decrease by $1$. On the other hand,
the node degree of a $\assocprodlambda$ product may decrease arbitrarily,
even if $j=\infty$:
this is a consequence of the non-homogeneity of the Baxter axiom~\eqref{E:rb},
and is not related to whether $\betaxy$ is quasi-idempotent or not.

In summary, we have
\begin{equation*}
\F_{\infty,j}(n_1,*)\assocprodlambda\F_{\infty,j}(n_2,*)\subseteq\F_{\infty,j}(n_1+n_2,*)
\end{equation*}
but only 
\begin{equation*}
\F_{2,j}(n_1,*)\assocprodlambda\F_{2,j}(n_2,*)\subseteq \F_{2,j}(n_1+n_2,*)\oplus \F_{2,j}(n_1+n_2-1,*)
\end{equation*}
and
\begin{equation*}
\Fxy(*,m_1)\assocprodlambda\Fxy(*,m_2)\subseteq\bigoplus_{\ell\leq m_1+m_2}\Fxy(*,\ell).
\end{equation*}
Thus,
the decomposition $\F_{\infty,j}=\bigoplus_{n\geq 0}\F_{\infty,j}(n,*)$ is an algebra grading and  the subspaces $\bigoplus_{\ell\leq m}\F_{2,j}(*,\ell)$ form an algebra filtration. In addition, 
the subspaces $\bigoplus_{\ell\leq m}\Fxy(*,\ell)$ form an algebra filtration, for
any $i,j$.
 
 The map $\betaxy:\Fxy\to\Fxy$ behaves as follows: 
 \begin{equation*}
 \beta_{i,\infty}\bigl(\F_{i,\infty}(n,m)\bigr)\subseteq \F_{i,\infty}(n,m+1)
 \quad\text{and}\quad
  \beta_{i,2}\bigl(\F_{i,2}(n,m)\bigr)\subseteq 
  \begin{cases}
\F_{i,2}(n,1) & \text{ if }m=0 \\
\F_{i,2}(n,m)         & \text{ if }m>0.
\end{cases}
\end{equation*}

\medskip

 For each $n\geq 1$ and $m\geq 0$, consider the dimensions
\begin{align*}
\fxy(n,m) &= \dim_\field\Fxy(n,m), &
\fxy(n,*) &= \dim_\field\Fxy(n,*), \\
\fxy(*,m) &= \dim_\field\Fxy(*,m), &
\fxy(k) &= \dim_\field\Fxy(k).
\end{align*}
Our main goal is to compute these dimensions  and to describe how they relate to each other as $i$ and $j$
vary over $\{2,\infty\}$.

Consider first the case $m=0$. A tree $t$ with $\degnode(t)=0$ has
only one internal node (the root), and this one is labeled by $0$. For
each $n\geq 1$ there is one such tree in $\T_{\infty,j}$, namely $\generator[n]$, with
bidegree $(n,0)$. Therefore, $\fii(n,0)=\fiw(n,0)=1$ for all $n\geq
1$. Similarly, 
\begin{equation*}
\fwi(n,0)=\fww(n,0)=\begin{cases}
  1 & \text{ if }n=1, \\
  0 & \text{ if }n>1.
\end{cases}
\end{equation*}

Unless explicitly stated, {}from now on we restrict our attention  to $n,m\geq 1$. 

Given
sequences $a(n)$ and $b(n,m)$ defined for $n,m\geq 1$,  the {\em binomial transforms} of
$a(n)$ and $b(n,m)$, respectively, are the sequences defined by
\begin{equation} \label{defBinomialTransform}
  \begin{aligned}
    \BT(a)(n) &= \sum_{k=1}^n \binom{n-1}{k-1}\, a(k), \\
    \BT^2(b)(n,m) &= \sum_{k=1}^n\sum_{\ell=1}^m \binom{n-1}{k-1}
    \binom{m-1}{\ell-1} \,b(k,\ell).
  \end{aligned}
\end{equation}
Clearly, the binomial transform $ \BT^2(b)$ can be computed as a double binomial transform, in any order: if we let $b^1_n(m) =\BT\bigl(b(\cdot,m)\bigr)(n)$ and $b^2_m(n) =\BT\bigl(b(n,\cdot)\bigr)(m)$, then
\begin{equation*}
  \BT^2(b)(n,m) = \BT(b^1_n)(m) = \BT(b^2_m)(n).
\end{equation*}

The following result says that as $i$ and $j$ vary, the dimensions of the homogeneous components $\Fxy(n,m)$ can be determined from the
dimensions of $\Fww(n,m)$ by applying binomial transforms.

\begin{proposition} \label{P:bin-tran}
 Consider the sequences $\fxy(n,m)$ for $n,m\ge 1$. We have 
   \begin{align*}
    \fwi(n,m) &= \BT\bigl(\fww(n,\cdot)\bigr)(m) \\
    \fiw(n,m) &= \BT\bigl(\fww(\cdot,m)\bigr)(n) \\
    \fii(n,m) &= \BT^2(\fww)(n,m) 
               = \BT\bigl( \fwi(\cdot, m) \bigr) (n)
               = \BT\bigl( \fiw(n, \cdot) \bigr) (m)
%    \fii(n,m) &= \BT(g_m)(n), \quad\text{where $g_m(\ell)=\BT\bigl(\fww(\ell,\cdot)\bigr)(m)$}
%    \fwi(n,m) &= \sum_{\ell=1}^m \binom{m-1}{\ell-1}\, \fww(n,\ell), \\
%    \fiw(n,m) &= \sum_{\ell=1}^n \binom{n-1}{\ell-1}\, \fww(\ell,m), \\
%    \fii(n,m) &= \sum_{\ell=1}^n\sum_{k=1}^m \binom{n-1}{\ell-1}
%    \binom{m-1}{k-1} \,\fww(\ell,k).
  \end{align*}
%  That is, it is possible to move around the spaces generated by the
%  diagram~\eqref{four-spaces} by computing binomial transforms of the
%  dimensions.
\end{proposition}

\begin{proof} Consider  the passage from $\fww$ to $\fiw$. Recall the map
$\tiw\onto\tww$ described in Proposition~\ref{P:can-maps}.
Let $t\in\tww(k,m)$. This is a tree with $k$ angles and $\degnode(t)=m$.
For each composition of $n$ with $k$ parts there is one tree $\hat{t}\in\tiw(n,m)$ in the fiber over $t$ of the map  (make the parts of the composition be the angle labels of $\hat{t}$). Since the number of such compositions is $\binom{n-1}{k-1}$, we obtain
\begin{equation*}
\fiw(n,m)=\sum_{k=1}^{n}\binom{n-1}{k-1}\fww(k,m)= \BT\bigl(\fww(\cdot,m)\bigr)(n).
\end{equation*}
The other cases are similar.
\end{proof}

The dimensions $\fww$ of the homogeneous components of
$\Fww$ admit very explicit descriptions, using the
bijections from Section~\ref{S:comb-trees}.

\begin{proposition} \label{P:fww}
  For $n\ge 1$ and $m\ge 0$, the dimension of the homogeneous
  components of $\Fww$ are given by:
  \begin{enumerate}[(i)]
  \item $\displaystyle{\fww(n,m) = \#\RP(n,m) = 
        C(m)\binom{m+1}{n-m}}$;
  \item $\fww(n,\ast) = 2\times \#\huMP(n)$;
  \item $\fww(\ast,m) = 2^{m+1}C(m)$. 
  \end{enumerate}
\end{proposition}

\begin{proof}
  Using Corollary~\ref{C:bijections} $(i)$, the equality
  $\fww(n,m)=\#\RP(n,m)$ is immediate. To count the number of these
  paths, remove the $n-m$ diagonal steps from one of those paths, the
  remaining steps can be assembled into an underdiagonal path from
  $(0,0)$ to $(m,m)$ with horizontal and vertical steps only. It is
  well-known that the number of such paths is the Catalan number
  $C(m)$~\cite[Exercise 6.19.h]{St99}. To reconstruct the given path
  from the Catalan path, since a diagonal step can only be followed by
  a horizontal step, there are $m+1$ possible places to distribute the
  $n-m$ diagonal steps: exactly before one of the horizontal steps, or
  in the last position.  This is $\binom{m+1}{n-m}$ possibilities, and
  thus the total number of paths is
  \begin{equation*}
    C(m)\binom{m+1}{n-m}
  \end{equation*}
  as claimed.
  
  Part $(ii)$ is a reformulation of Corollary~\ref{C:bijections}
  $(iii)$. For part $(iii)$ write, using $(i)$,
  \begin{equation*}
    \fww(*,m)=\sum_{n\ge 1}\fww(n,m)=\sum_{n\ge 1}
    C(m)\binom{m+1}{n-m}=C(m)2^{m+1},
  \end{equation*}
  by the binomial theorem.
\end{proof}
Observe that $\fww(n,m)$ is non-zero only in the region $m\le n\le
2m+1$.

%%% use only one kind of paths (always motzkin). Didn't understand "this justifies..." Also, is it exhaustive?

The dimensions $\fiw$ of the homogeneous components of $\Fiw$ also
admit simple combinatorial descriptions, in addition to the
descriptions in terms of decorated trees or in terms of the binomial
transform (Proposition~\ref{P:bin-tran}).

\begin{proposition} \label{P:fiw}
  For $1\le m\le n$, the dimension of the homogeneous components of
  $\Fiw$ are given by:
  \begin{enumerate}[(i)]
  \item $\displaystyle{\fiw(n,m)= \#\SP(n,m) = C(m)\binom{n+m}{n-m}}$;
  \item $\fiw(n,*) = 2\times\#\PT(n)$, which are the large
    Schr\"oder numbers;   % twice the number of planar rooted trees with $n+1$ leaves. This is
    % sequence A006318.
  \item $\fiw(*,m)$ is infinite.
  \end{enumerate}
\end{proposition}

\begin{proof}
  Use Corollary~\ref{C:bijections} $(i)$ to conclude
  $\fiw(n,m)=\#\SP(n,m)$. We proceed similarly as before to count this
  number. Given a path in $\SP(n,m)$, after removing the $n-m$
  diagonal steps, we get a Catalan path. To reconstruct the initial
  path we need to distribute the $n-m$ diagonal steps in $2m+1$
  places: before one of the $m$ horizontal steps, before one of the
  $m$ vertical steps, or in the last position of the path. Since there
  can be many consecutive diagonal steps in each place, the total
  number is
  \begin{equation*}
    C(m)\binom{(2m+1)+(n-m)-1}{n-m} = C(m)\binom{n+m}{n-m}.
  \end{equation*}
  Part $(ii)$ is again a direct consequence of
  Corollary~\ref{C:bijections} $(ii)$. Part $(iii)$ is clear since the
  decorations in the angles are arbitrary.
\end{proof}
%\begin{proof}
%  A bijection for the part $(i)$ is easy. Given a tree $t\in\tiw(n,m)$,
%  substitute the decorations $j$ in each angle for $j-1$ intermediate
%  leaves in the corresponding node. The resulting tree will have $n$
%  angles, hence $n+1$ leaves and $m$ or $m+1$ nodes, according to
%  wether the label in the root is $1$ or $0$. For example,
%  \begin{equation*}
%\vc{\rbtree{A}{levelsep=3ex}{r}{\rootlbl{\;1}}{%
%  \rbleaf{B}%
%  \rbtree{C}{}{r}{}{%
%    \rbleaf{D}%
%    \rbleaf{E}%
%  }%
%  \rbleaf{F}%
%}%
%\angled{B}{C}{1}%
%\angled{C}{F}{2}%
%\angled{D}{E}{3}%
%} \qquad\longmapsto\qquad
%\vc{%
%  \pstree[levelsep=3ex,thistreesep=8pt]%
%    {\TC*}%
%    {\TC* {\pstree[levelsep=4ex,thistreesep=6pt]
%        {\TC*} {\TC* \TC* \TC* \TC*}}
%      \TC*
%      \TC*}}
%  \end{equation*}
%For the reverse process
%  use the number of nodes to decide the label in the root. The number
%  of these trees is easily computed via the binomial transform of
%  $\fww$. The expression in term of lattice paths comes from the
%  application of algorithms~\ref{treetopath} and~\ref{pathtotree} to
%  the rooted planar trees.

%  For part $(ii)$ use the same substitution as in part $(i)$, the
%  label $0$ or $1$ in the root gives two copies of the rooted planar
%  trees with $n+1$ leaves.

%  Part $(iii)$ is clear, since the decorations in the angles are arbitrary.
%\end{proof}

\begin{remark}
  The three previous propositions show that the small Schr\"oder
  numbers are the binomial transform of the $(h,u)$-colored
  Motzkin numbers, a result stated by D. Callan in~\cite{Sl}.
\end{remark}

Among the sequences $\fxy(k)$, the case of $\fiw$ again proves  to be
interesting combinatorially.

\begin{proposition} \label{P:fxw(k)}
  For $k\ge 1$, the dimensions of the subspaces $\Fww(k)$ and
  $\Fiw(k)$ are
  \begin{enumerate}[(i)]
  \item $\fww(k) = \#\RMP(k)$,
  \item $\fiw(k) = \#\MP(k)$, which are the Motzkin numbers.
  \end{enumerate}
\end{proposition}

\begin{proof}
  This is a restatement of Proposition~\ref{P:fiw-k}.
\end{proof}

  The dimensions $\fwi(n,m)$ and $\fii(n,m)$ do not seem to admit any simpler description than as the iterated
  binomial transforms of $\fww$. We mention that the sequences $\fwi(*,m)$ and
  $\fwi(k)$
  appear in~\cite{Sl} as A082298 and A025243, respectively, while the
  dimensions $\fwi(n,*)$, $\fii(n,*)$ and $\fii(*,m)$ are 
  infinite.
  
  Some of the sequences from Propositions~\ref{P:fww}, \ref{P:fiw},
  and~\ref{P:fxw(k)} also appear in~\cite{Sl} as $\fww(n,m)$: A068763;
  $\fww(n,\ast)$: A071356; $\fww(\ast,m)$: A025225;
  $\fiw(n,\ast)$: A006318; $\fww(k)$: A007477; $\fwi(k)$: A025243.
  
 Table~\ref{T:dimensions} summarizes the results of this section. 

\begin{table}[!ht]
\centering
\renewcommand\arraystretch{1.5}
\begin{tabular}{|>{$}c<{$}||>{$}c<{$}|>{$}c<{$}|>{$}c<{$}|>{$}c<{$}|}
\hline
\rule{0pt}{14pt}
\rule[-7pt]{0pt}{0pt}
 & \mathord{\cdot}(n,m)
 & \mathord{\cdot}(n,*) 
 & \mathord{\cdot}(*,m) 
 & \mathord{\cdot}(k) 
\\
\hhline{|=::====|}
\rule{0pt}{18pt}
\rule[-12pt]{0pt}{0pt}
\parbox[c]{70pt}{\centering $\fww$ \\ 
   \tiny $(\lfloor n/2 \rfloor\le m\le n)$} 
 & 
 \displaystyle{C(m)\binom{m+1}{n-m}}
 & 
\vborder{4pt}{\parbox[c]{120pt}{\small\centering twice the number
      of \\$(h,u)$-colored Motzkin  paths}}
 & 
 2^{m+1}C(m)
 & 
 \#\RMP(k)
\\
\hhline{|-||----|}
\vborder{4pt}{\parbox[c]{70pt}{\centering $\fiw$ \\ 
  \tiny $(0 \le m\le n)$}}
 & 
 \vborder{4pt}{\ensuremath{\displaystyle{C(m)\binom{n+m}{n-m}}}}
 & 
 \parbox[c]{80pt}{\small\centering large
   Schr\"oder\\ number} 
 & 
 \infty
 &
 \parbox[c]{50pt}{\small\centering Motzkin \\ number}\\
%\hhline{|-||----|}
%%
%\rule{0pt}{18pt}
%\rule[-12pt]{0pt}{0pt}
%\parbox[c]{70pt}{\centering $\fwi$ \\ 
%\tiny $(1 \le n\le 2m+1)$}
%& \BT\bigl(\fww(n,\cdot)(m)\bigr) & \infty & \text{A082298}  &
%\text{A025243} \\
%\hhline{|-||----|}
%%%
%\rule{0pt}{18pt}
%\rule[-12pt]{0pt}{0pt}
%\parbox[c]{70pt}{\centering $\fii$ \\ 
%\tiny $(n\ge 1;\; m\ge 0)$}
%& \BT^2(\fii)(n,m) & \infty & \infty & \parbox[c]{80pt}{\small\centering
%  generating\\
%function~\eqref{fiik}} \\
\hline
\end{tabular}
\caption{Dimensions of the free \rba s with a quasi-idempotent operator} \label{T:dimensions}
\end{table}

\subsection{Generating series of the free \rba s}\label{S:generating}

Given sequences $a(n)$ and $b(n,m)$ defined for $n,m\geq 1$, consider their  generating
  functions $A\in\field\llbracket x\rrbracket$ and
$B\in\field\llbracket x,y\rrbracket$, defined  by
\begin{equation*}
A(x)=\sum_{n\ge 1}a(n)\,x^n,\quad  
B(x,y) = \sum_{n\ge 1}\sum_{m\ge 1} b(n,m)\, x^n y^m.
\end{equation*}
The binomial transform~\eqref{defBinomialTransform} has a simple expression in terms of generating
functions. If $c = \BT(a)$, then the generating functions $C$ and $A$
are related by
\begin{equation} \label{gen-BT}
  C(x) = A\Bigl( \frac{x}{1-x} \Bigr).
\end{equation}

Below we find a closed expression for the generating function
of the sequence $\fww(n,m)$, and then we use 
Proposition~\ref{P:bin-tran} to obtain the generating
functions for the sequences $\fwi(n,m)$,  $\fiw(n,m)$, and  $\fii(n,m)$.

\begin{proposition}\label{P:gen-fij}
  The generating functions for the sequences $\fxy(n,m)$, with
  $i,j\in\{2,\infty\}$, are as follows:
  \begin{align*}
    \Gww(x,y) &= (1+x)\,f\bigl(xy(1+x)\bigr), \\
    \Gwi(x,y) &= (1+x)\, f\Bigl( \frac{xy(1+x)}{1-y} \Bigr), \\
    \Giw(x,y) &= \frac{1}{1-x} \,f\Bigl( \frac{xy}{(1-x)^2} \Bigr), \\
    \Gii(x,y) &= \frac{1}{1-x}\, f\Bigl( \frac{xy}{(1-x)^2(1-y)} \Bigr),
  \end{align*}
  where $f(u)=\sum_{n\ge 1}C(n)u^n=(1-\sqrt{1-4u}\,) / (2u) - 1$ is
  the generating function for the Catalan numbers.
\end{proposition}

\begin{proof}
  The  last three formulas follow from the first in view of~\eqref{gen-BT}. To verify the first we compute
    \begin{align*}
    \Gww(x,y) &=\sum_{m\ge 1}\sum_{n=m}^{2m+1} C(m)\binom{m+1}{n-m}x^ny^m \\
    &= \sum_{m\ge 1} C(m)(xy)^m \sum_{k=0}^{m+1} \binom{m+1}{k} x^k \\
    &= \sum_{m\ge 1} C(m)(xy)^m(1+x)^{m+1} \\
    &= (1+x)f\bigl(xy(1+x)\bigr),
  \end{align*}
  as claimed.
\end{proof}

The generating functions for the sequences $\fxy(k)$ are obtained by setting $x=y$ in Proposition~\ref{P:gen-fij}.

\section{Free algebra with an idempotent morphism}\label{S:idem}

Let $A$ be an algebra and $\beta:A\to A$ an idempotent morphism of algebras:
\begin{equation*}
\beta(a)\beta(b)=\beta(ab), \qquad \beta\bigl(\beta(a)\bigr)=\beta(a).
\end{equation*}
Then also
\begin{equation*}
\beta(a)\beta(b)=\beta\bigl(\beta(a)b+a\beta(b)-ab\bigr),
\end{equation*}
so $\beta$ is an idempotent \rb\ operator of weight $\lambda=-1$. 
We may thus
consider the full subcategory $\Mori$ of $\Riw^{-1}$ whose objects are triples 
$(A,x,\beta)$
where $A$ is an algebra, $x\in A$, and $\beta:A\to A$ is an
idempotent morphism of algebras. We refer to the initial object in this category as
the free algebra on one generator with an idempotent morphism. 

Similarly, by the  free algebra with an idempotent morphism and an idempotent generator we mean the initial object in the full subcategory $\Moriw$ of $\Mori$
whose objects satisfy $x^2=x$.

Since $\Mori$ is a subcategory of $\Riw^{-1}$,
there is a unique morphism of \rba s from
$(\Fiw^{-1},\generator,\betaiw)$, the initial object in the category $\Riw^{-1}$, to 
the initial object in $\Mori$. Similarly, there is a unique morphism of \rba s from
$(\Fww^{-1},\generator,\betaww)$ to the initial object in $\Moriw$.
We proceed to construct these initial objects and to describe 
these canonical morphisms  in explicit terms.

\subsection{Construction of the free algebras with an idempotent morphism}\label{S:const-idem}

Let $\Mi=\field\langle x_0,x_1\rangle$ be the free associative
algebra in two variables $x_0$ and $ x_1$. For consistence with the preceding constructions, we stick to the world of non-unital algebras (thus, we leave the constant polynomials out). As explained in Remark~\ref{R:unital} (below), this is not
an essential restriction in this context.

Let $\betai:\Mi\to \Mi$ be the
unique morphism of algebras  such that
\begin{equation*} 
  \betai(x_0) = x_1,\quad \betai(x_1) = x_1.
\end{equation*}

\begin{proposition}
The initial object in the category
$\Mori$ is $(\Mi,x_0,\betai)$. 
\end{proposition}
\begin{proof}
 The map $\betai$  is 
  idempotent on the generators, hence everywhere. Therefore, the
  object $(\Mi,x_0,\betai)$ belongs to the category $\Mori$.

  Let $(A,x,\beta)$ be another object  in $\Mori$. 
  Let  $\varphi:\Mi\to A$ be the unique morphism of algebras
  such that
  \begin{equation*}
    \varphi(x_0) = x,\quad \varphi(x_1) = \beta(x).
  \end{equation*}
  We have
  \begin{equation*}
    \beta\varphi(x_0)=\beta(x)=\varphi(x_1)=\varphi\betai(x_0)
  \end{equation*}
  and
  \begin{equation*}
    \beta\varphi(x_1)=\beta\beta(x)=\beta(x)=\varphi(x_1)=\varphi\betai(x_1),
  \end{equation*}
  since $\beta$ is idempotent. Since all these maps are morphisms of algebras,
  we get that
  \begin{equation*}
    \beta\varphi=\varphi\betai.
  \end{equation*}
  Thus $(\Mi,x_0,\betai)$ is
  the initial object in $\Mori$.
  \end{proof}

  In particular, $(\Mi,\betai)$ is a \rba. Let $\pii:\Fiw^{-1}\to\Mi$ be the unique morphism of \rba s such that
   \begin{equation}\label{E:gen-0}
    \pii\bigl(\,\generator\,\bigr) = x_0.
  \end{equation}
  We have
   \begin{equation} \label{E:gen-1}
  \pii(\betagenerator)=\pii\Bigl(\betaiw\bigl(\generator\bigr) \Bigr) =
  \betai\Bigl(\pii\bigl(\generator\bigr) \Bigr) = \betai(x_0) = x_1.
  \end{equation}
  
  More generally:
  \begin{lemma}\label{L:pi} Let $t\in\tiw$, $t_1,\ldots,t_{n}\in\widehat{\T}_{\infty,2}^+$ the subtrees of $t$ rooted at the children of the root, and $i_1,\ldots,i_{n-1}$ the labels of the angles between these children, as in~\eqref{E:def-H}. 
   Then
  \begin{equation}\label{E:pi}
\pii(t)=\begin{cases}
x_1^{\degangle(t_1)}
  x_0^{i_1} \cdots x_0^{i_{n-1}} x_1^{\degangle(t_{n})} & \text{ if  $t\in\tiw^0$, } \\
 x_1^{\degangle(t)}         & \text{ if $t\in\tiw^+$.}
\end{cases}
\end{equation}
\end{lemma}
\begin{proof}  
We argue by induction on the bidegree of $t$, starting from~\eqref{E:gen-0}.
According to~\eqref{E:can-exp}, we have 
\begin{equation*}
   \pii(t) = \betai^a\bigl( \pii(t_1)  x_0^{i_1} 
                 \pii(t_2) x_0^{i_2}  \cdots
                x_0^{i_{n-1}}  \pii(t_n)\bigr),
  \end{equation*}        
  where $a=0$ if  $t\in\tiw^0$ and $a=1$ if  $t\in\tiw^+$.        
   Now, the trees $t_{k}$ belong to $\tiw^+$ and have smaller degree than $t$, so  by induction hypothesis $\pii(t_k)= x_1^{\degangle(t_k)}$ (if $t_1=\point$ or $t_n=\point$, then they do not appear in the above expression).
 Substituting above we get
 \begin{equation*}
   \pii(t)= \betai^a\bigl( x_1^{\degangle(t_1)}  x_0^{i_1} 
   x_2^{\degangle(t_2)} x_0^{i_2}  \cdots
   x_0^{i_{n-1}}  x_n^{\degangle(t_n)}\bigr).
 \end{equation*}
Using that $\betai$ is a morphism and $\betai(x_0)=x_1$ we obtain~\eqref{E:pi}.                 
\end{proof}

In order to describe the kernel of $\pii$, we introduce the following relation among decorated trees. Recall that the elements of $\tiw$ are trees whose root label is $0$ or $1$ and the only other decorations are in the angles.
%%% Change notation below to agree with new grafting notation.
%%% Also, perhaps say t_i\sim s_i instead of \deg(t_i)=\deg(s_i)
Given $t$ and $s$ in $\tiw$,  write $t\sim s$ if
the following conditions hold:
\begin{enumerate}[($a$)]
\item $t$ and $s$ have the same root label.
\item If the root label is $1$, then 
\begin{enumerate}[($i$)]
\item $\degangle(t)=\degangle(s)$.
\end{enumerate}
If the root label is $0$, then 
\begin{enumerate}[($i$)]
\item $t$ and $s$ have  the same number  of children of the root, 
\item $\degangle(t_k)=\degangle(s_k)$ for all $k=1,\ldots,n$,
\item $i_k=j_k$ for all $k=1,\ldots,n-1$,
\end{enumerate}
  where $t_1,\ldots,t_n$ are  the subtrees of $t$ rooted at the children of the root, and $i_1,\ldots,i_{n-1}$ are the labels of the angles between these children, as in~\eqref{E:def-H}, and similarly for $s_k$, $j_k$, and $s$.
\end{enumerate}
For example, the following trees are related:
\begin{equation*}
  \tree{\children{A}{,xbbh=\lblheight,xbbr=.4\treesepone,%
                     xbbd=\levelsep,xbbl=.4\treesepone}{2}{1.4\treesepone}{\levelsep}%
    \rput[t](AChild1){%
      \children{B}{}{2}{.8\treesepone}{\levelsep}}%
    \rput[t](AChild2){%
      \children{C}{}{2}{.8\treesepone}{\levelsep}}%
    \lbl[bl]{ARoot}{0pt}{\lblsep}{0}%
    \angled{AChild1}{AChild2}{2}%
    \angled{BChild1}{BChild2}{3}%
    \angled{CChild1}{CChild2}{1}%
  }
\quad \sim \quad
  \tree{\children{A}{,xbbh=\lblheight,xbbr=.4\treesepone,%
                     xbbd=2\levelsep,xbbl=.4\treesepone}{2}{1.4\treesepone}{\levelsep}%
    \rput[t](AChild1){%
      \children{B}{}{2}{.8\treesepone}{\levelsep}}%
    \rput[t](AChild2){%
      \children{C}{}{2}{.8\treesepone}{\levelsep}}%
    \rput[t](BChild2){%
      \children{D}{}{2}{.8\treesepone}{\levelsep}}%
    \lbl[bl]{ARoot}{0pt}{\lblsep}{0}%
    \angled{AChild1}{AChild2}{2}%
    \angled{BChild1}{BChild2}{1}%
    \angled{CChild1}{CChild2}{1}%
    \angled{DChild1}{DChild2}{2}%
  }.
\end{equation*}

\begin{proposition}\label{P:Mi} The vector subspace $I$ of $\Fiw^{-1}$ spanned by the differences $t-s$ for $t\sim s$ in $\tiw$ is a \rb\ ideal, and
  \begin{equation*}
    \Mi\cong \Fiw^{-1}/I.
  \end{equation*}
\end{proposition}

\begin{proof} Consider the morphism $\pii:\Fiw^{-1}\to\Mi$.
According to~\eqref{E:gen-0} and~\eqref{E:gen-1},
 the algebra generators of $\Mi$ are in its image, so this map  is surjective.
 Now, by Lemma~\ref{L:pi}, $\pii$ sends the basis $\tiw$ of $\Fiw^{-1}$ to the basis of words in $x_0$ and $x_1$ of $\Mi$, from which it follows that the kernel of $\pii$ is the subspace spanned by the differences $t-s$, for all $t,s\in\tiw$ with $\pii(t)=\pii(s)$. Now, the latter occurs precisely when $t\sim s$, by~\eqref{E:pi}.
\end{proof}

The free algebra with an idempotent morphism and an idempotent
generator admits a similar description. Let $ \Mw$ be the following quotient of
$\Mi$:
  \begin{equation*}
    \Mw = \frac{\field\langle x_0,x_1\rangle}{\langle x_0^2=x_0,x_1^2=x_1\rangle}.
  \end{equation*}
The morphism $\betai:\Mi\to\Mi$ descends to the quotient, giving rise to another idempotent
morphism $\betaw:\Mw\to\Mw$, which turns $(\Mw,\betaw)$ into a \rba\ of weight $\lambda=-1$. Since $x_0$ is now idempotent, $(\Mw,x_0,\betaw)$ is an object of the category $\Rww^{-1}$, and there is a unique
morphism of \rba s  $\piw:\Fww^{-1}\to\Mw$ such that
$\piw(\generator)=x_0$. These  fit into a commutative diagram of 
of morphisms of \rba s:
\begin{equation} \label{idem-mor}
  \begin{gathered}
    {\psset{nodesep=2pt}%
      \begin{psmatrix}[rowsep=20pt,colsep=30pt]
        \Fiw^{-1} & \Mi \\ 
        \Fww^{-1} & \Mw
        \ncline{->>}{1,1}{1,2}\naput{\scriptstyle{\pii}}
        \ncline{->>}{1,1}{2,1}
        \ncline{->>}{1,2}{2,2}
        \ncline{->>}{2,1}{2,2}\naput{\scriptstyle{\piw}}
      \end{psmatrix}%
    }%
  \end{gathered}
\end{equation}

  \begin{remark}\label{R:unital} For general \rba s, the distinction between the free algebras in the unital case and in the non-unital case is important. In this paper we have dealt with non-unital algebras only. On the other hand, this becomes a minor point when considering the free algebras with an idempotent morphism,
as  the preceding constructions show. Namely, to obtain the free unital algebras
with an idempotent morphism, simply throw in the constant polynomials to the
spaces $\Mi$ and $\Mw$, and extend the morphisms $\betai$ and $\betaw$  so that they preserve the unit element.
  (Note that for an arbitrary \rb\ operator such an extension does not produce
another \rb\ operator.)
\end{remark}

  \subsection{Bigrading and dimensions of the free algebras with an idempotent morphism}

 Consider the  bigrading of the polynomial algebra $\Mi$  defined as follows. Any monomial $\mu$ in $x_0$ and $x_1$ can be uniquely written as $\mu=x_0^{i_0}x_1^{j_1}x_0^{i_1}x_1^{j_2}\ldots x_1^{j_k}x_0^{i_k}$ with $i_0,i_k\geq 0$ and all other exponents  $i_h,j_h>0$. Then set 
  \begin{equation}\label{E:deg-M}
\deg(\mu)=(i_0+j_1+i_1+j_2+\cdots+j_k+i_k,\,k).
\end{equation}
In particular,
\begin{equation*}
\deg(x_0)=(1,0) \text{ \ and \ } \deg(x_1)=(1,1).
\end{equation*}
This can be understood as follows:
 $\deg(\mu)=(n,m)$ if when writing $\mu\in\Mi$ as a word in  $x_0$ and $\betai(x_0)$, the symbol $x_0$ occurs exactly $n$ times and the symbol $\betai$ occurs {\em at least} $m$ times.

Let $\Mi(n,m)$, $\Mi(n,*)$, and $\Mi(*,m)$ be the corresponding homogeneous components. In analogy with the situation encountered for the algebras $\Fxy^{\lambda}$ in Section~\ref{S:comb-dim}, we have
 that the decomposition $\Mi=\bigoplus_{n\geq 0}\Mi(n,*)$ is an algebra grading, but the decomposition $\Mi=\bigoplus_{m\geq 0}\Mi(*,m)$ is not:
$x_1^2$ has degree $(2,1)$.
 On the other hand,  the subspaces $\bigoplus_{\ell\leq m}\Mi(*,\ell)$ form an algebra filtration.
  
The subspaces  $\bigoplus_{\ell\leq n}\Mw(\ell,*)$
define an algebra filtration on $\Mw$ (now the exponents $i_h$ and $j_h$ are at most $1$), and the quotient map $\Mi\to\Mw$ is filtration-preserving. In addition, the morphisms  
 $\pii:\Fiw^{-1}\to\Mi$ and $\piw:\Fww^{-1}\to\Mw$ are degree-preserving, in view of~\eqref{E:pi}, so all maps in~\eqref{idem-mor} are degree-preserving morphisms of algebras. We use the same notation for the various homogeneous components of $\Mw$, and  we obtain two algebra filtrations on $\Mw$. 
 
The map $\beta_i:M_i\to M_i$ behaves as follows:
 \[
  \betai\bigl(\Mi(n,m)\bigr)\subseteq \Mi(n,m+1)
  \quad\text{and}\quad
  \beta_2\bigl(M_2(n,m)\bigr)\subseteq
  \begin{cases}
M_2(n,1) & \text{ if }m=0 \\
 M_2(n,m)        & \text{ if }m>0.
\end{cases}
  \]

 Let $\mi(n,m)=\dim_\field \Mi(n,m)$
 and $\mw(n,m)=\dim_\field \Mw(n,m)$
  be the dimensions of the homogeneous components of bidegree $(n,m)$.
  It is easy to see that, for any $n,m\geq 0$,
   \begin{equation} \label{freealgidempot}
 \mi(n,m)=\binom{n+1}{2m} \text{ \ and \ }  \mw(n,m)= 
    \begin{cases}
      2 &\text{if $n=2m$ and $(n,m)\neq(0,0)$,} \\
      1 &\text{if $|n-2m|=1$ or $(n,m)=(0,0)$,}\\
      0 &\text{otherwise.}
    \end{cases}
\end{equation}
 Note that, in analogy with the situation encountered for free \rba s in
Proposition~\ref{P:bin-tran}, these dimensions are
   related by a binomial transform:
  \begin{equation} \label{E:bin-tran-mor}
    \mi(n,m) = \BT\bigl(\mw(\cdot,m)\bigr)(n)=\sum_{i=1}^n\binom{n-1}{i-1}\mw(i,m).
  \end{equation}
  This assertion boils down to Pascal's identity for binomial coefficients.

%Below, we take $\deg(x_1)=(0,1)$ which leads to nice formulas but does not
%agree with the degree defined on $\Fiw$.

%  Consider the usual grading of the polynomial algebras $\Mi$ and $\Mw$, for which the variables $x_0$ and $x_1$ have
% bidegree $(1,0)$ and $(0,1)$, respectively.  Let $\mi(n,m)=\dim_\field \Mi(n,m)$
% and $\mw(n,m)=\dim_\field \Mw(n,m)$
%  be the dimensions of the homogeneous components of bidegree $(n,m)$.
%  We have
%   \begin{equation} \label{freealgidempot}
% \mi(n,m)=\binom{n+m}{n} \text{ \ and \ }  \mw(n,m)= 
%    \begin{cases}
%      2 &\text{if $n=m$,} \\
%      1 &\text{if $|n-m|=1$,}\\
%      0 &\text{otherwise.}
%    \end{cases}
%\end{equation}
%The first formula is clear; the second follows from the fact that the set of words in alternating letters $x_0$ and $x_1$ is a linear basis
%of $\Mw$. Note that, in analogy with the situation encountered for free \rba s in
%Proposition~\ref{P:bin-tran}, these dimensions are
%   related by a double binomial transform:
%  \begin{equation} \label{E:bin-tran-mor}
%    \mi(n,m) = \BT^2(\mw)(n,m).
%  \end{equation}
%  This assertion boils down to Vandermonde's identity for binomial coefficients.

\subsection{Generating series of the free algebras with an idempotent morphism}
 We consider the unital version of these algebras. As explained in Remark~\ref{R:unital}, this simply amounts to adding one to the generating series of the non-unital versions.
  The generating functions for the sequences $\mi(n,m)$ and $\mw(n,m)$, $n,m\geq 0$, are easily seen to be
  \begin{equation*}
    \mathsf{M}_\infty(x,y) = \frac{1-x+xy}{(1-x)^2-x^2y}    \text{ \ and \ }  \mathsf{M}_2(x,y) = \frac{(1+x)(1+xy)}{1-x^2y}.
  \end{equation*}
Note that ${\displaystyle \mathsf{M}_\infty(x,y) = \mathsf{M}_2\Bigl( \frac{x}{1-x},y\Bigr)}$, in agreement with~\eqref{gen-BT} and~\eqref{E:bin-tran-mor}.

\section{Connections with dendriform trialgebras and dialgebras}\label{S:dendri}

\subsection{The free dendriform trialgebra and the free dendriform dialgebra}\label{S:dend-free}

Dendriform dialgebras and trialgebras were introduced by Loday~\cite{L} and
Loday and Ronco~\cite{LR}. For our purposes it is convenient to consider the following notion.

\begin{definition}\label{D:dend-tri} Fix $\lambda\in\field$, a
{\em $\lambda$-dendriform trialgebra} $D$ is a vector space with three
binary operations $\prec$, $\succ$, and $\cdot$, verifying for all $x,y,z\in D$,
\smallskip
\begin{gather*}
\renewcommand\minalignsep{20pt}
\begin{aligned}
  (x\ast y)\succ z &= x\succ(y\succ z), %
    & (x\succ y)\cdot z &= x\succ (y\cdot z), \\
     (x\succ y)\prec z &= x\succ (y\prec z), %
  &  (x\prec y)\cdot z &= x\cdot (y\succ z), \\
       (x\prec y)\prec z &= x\prec(y\ast z), %
      & (x\cdot y)\prec z &= x\cdot (y\prec z),\\
\end{aligned} \\
  (x\cdot y)\cdot z = x\cdot (y\cdot z),
\end{gather*}
where
\begin{equation*}
x\ast y = x\prec y + x\succ y + \lambda
(x\cdot y).
\end{equation*}
\end{definition}

For $\lambda=1$ we obtain the usual notion of dendriform trialgebras.
For any $\lambda\in\field$, any $\lambda$-dendriform trialgebra may be turned into a $1$-dendriform trialgebra by means of the transformation
\begin{equation*}
(D,\prec, \succ,\cdot) \mapsto (D,\prec,\succ, \lambda\mathord{\cdot})
\end{equation*}
(multiplying the last operation by $\lambda$).
If $\lambda\neq 0$, this transformation is invertible, but if $\lambda=0$
a truly distinct notion arises. This notion is closely related to, but not the same as, that of dendriform dialgebras~\cite{L}. Specifically, any $0$-dendriform trialgebra may be turned into a dendriform dialgebra by means of the transformation
\begin{equation*}
(D,\prec, \succ,\cdot) \mapsto (D,\prec,\succ)
\end{equation*}
(forgetting the last operation). This transformation is not invertible.

Let $\DT^\lambda$ denote the category whose objects are pairs $(D,x)$ where $D$ is a $\lambda$-dendriform trialgebra and $x\in D$, and whose morphisms are maps that preserve the operations and the distinguished elements.

The initial object in $\DT^1$ (the free dendriform trialgebra on one generator)
was constructed in explicit combinatorial terms in~\cite{LR}. A slight variant
of this construction leads to the initial object in $\DT^\lambda$ for any $\lambda\in\field$. 

Recall that $\PT$ denotes the set of rooted planar trees, and $\PT(n,m)$ consists of trees with $n+1$ leaves and $m$ internal nodes (Section~\ref{S:comb-trees}).

\begin{proposition}[\cite{LR}] \label{P:free-dend}
Let $\FDT$ be the vector space with basis consisting of the set $\bigoplus_{n,m\geq 1}\PT(n,m)$. Fix $\lambda\in\field$ and define operations on this space by means of the following recursions:
\begin{align}\label{E:prec}
x\prec y & = \Gr(x_1,\ldots,x_n\ast y),\\
\label{E:succ}
x\succ y & =\Gr(x\ast y_1,y_2,\ldots,y_m),\\
\label{E:cdot}
x\cdot y & =\Gr(x_1,\ldots,x_n\ast y_1,y_2,\ldots,y_m),\\
x\ast y & =x\prec y + x\succ y +\lambda (x\cdot y).
\end{align}
Let $\FDT^\lambda$ denote the space $\FDT$ endowed with the operations $\prec$, $\succ$, and $\cdot$. Then $(\FDT^\lambda,\gendend)$ is the initial object in $\DT^\lambda$.
\end{proposition}

In the above definitions, we have set $\H(x)=(x_1,\ldots,x_n)$ and
$\H(y)=(y_1,\ldots,y_m)$, and $\Gr$ and $\H$ stand for grafting  and
de-grafting of rooted planar binary trees: $\Gr$ is defined as
in~\eqref{E:def-G-many} and $\H$ is defined as in the first case of~\eqref{E:def-H}, ignoring all
labels in both definitions. In contrast to the grafting in~\eqref{E:def-G-many},
no normalization is required, since internal leaves are now allowed. The operation $\ast$ is defined on the larger space spanned by $\bigoplus_{n,m\geq 0}\PT(n,m)$ and the recursion starts with $x\ast\point=\point\ast x= x$.

\smallskip

Let us consider the analogous notions for dendriform dialgebras.
Let $\DD$ denote the category whose objects are pairs $(D,x)$ where $D$ is a dendriform dialgebra and $x\in D$, and whose morphisms are maps that preserve the operations and the distinguised elements.
The initial object in $\DD$ (the free dendriform dialgebra on one generator)
was constructed  in~\cite{L}. On the vector space
$\FDD$ with basis consisting of the set of rooted planar binary trees, two operations $\prec$ and $\succ$ are defined by means of formulas similar
to those in Proposition~\ref{P:free-dend}. The result $(\FDD,\gendend)$ is the initial object in $\DD$.

\subsection{Embedding dendriform trialgebras and dialgebras in \rba s}\label{S:embedding}

The following observation relates dendriform trialgebras and dialgebras to \rba s.

\begin{proposition}[\cite{poisson,E}]\label{P:dendri-baxter}
Let $(A,\beta)$ be a \rba\ of weight $\lambda$. Defining
\begin{equation*}
x\succ y=\beta(x)y,\qquad  x\prec y=x\beta(y), \text{ \ and \ }
x\cdot y= xy
\end{equation*}
one obtains a $\lambda$-dendriform trialgebra structure on $A$.
\end{proposition}

In view of Proposition~\ref{P:dendri-baxter}, we may turn the free \rba\ $\Fii^\lambda$ into a $\lambda$-dendriform trialgebra. Therefore, there is a unique
morphism of dendriform trialgebras
\begin{equation*}
\FDT^\lambda \to \Fii^\lambda
\end{equation*}
that sends $\gendend$ to $\generator$.
Ebrahimi-Fard and Guo used their construction of the free \rba\ to make the interesting observation that this map is injective~\cite{EG}. Below we derive the stronger
fact that the composite
\begin{equation*}
i:\FDT^\lambda \to \Fii^\lambda\onto \Fiw^\lambda
\end{equation*}
is still injective, and describe these map in explicit combinatorial terms.

\begin{proposition}\label{P:trialg-inj}
  The canonical morphism of dendriform trialgebras
  \begin{equation*}
    i:\FDT^\lambda\to\Fiw^\lambda,
  \end{equation*}
  sends any rooted planar tree $x\in\PT$ to the decorated tree
  $(f^0)^{-1}(x)\in\T^0_{\infty,2}$, where $f^0$ is the bijection of
  Proposition~\ref{P:bijections}.  In particular, $i$ is injective.
\end{proposition}

\begin{proof}
  Let $i':\FDT^\lambda\to\Fiw^\lambda$ by the map defined by
  $i'(x)=(f^0)^{-1}(x)$.  Then $i'(\gendend)=\generator$. We show
  below that $i'$ is a morphism of dendriform trialgebras; then, by
  uniqueness, $i'=i$.
  
  We proceed by induction on the number of nodes of $x$ and $y$,
  proving that $i'$ preserves the three operations on $\FDT^\lambda$.
  The equality $i'(x\succ x) = i'(x)\succ i'(x)$, where $x=\gendend$
  is immediate:
  \begin{align*}
    i'(\gendend\succ\gendend) &= i'\Bigl( 
    \tree{\children{A}{,xbbl=.4\treesepone,xbbd=.7\levelsep}{2}{.8\treesepone}{.7\levelsep}%
      \rput[t](AChild1){\children{B}{}{2}{.8\treesepone}{.7\levelsep}}%
%      \lbl[bl]{ARoot}{0pt}{\lblsep}{0}%
    } \Bigr) =
    \tree{\children{A}{,xbbh=\lblheight,xbbl=.4\treesepone,xbbd=.7\levelsep}{2}{.8\treesepone}{.7\levelsep}%
      \rput[t](AChild1){\children{B}{}{2}{.8\treesepone}{.7\levelsep}}%
      \lbl[bl]{ARoot}{0pt}{\lblsep}{0}%
      \angled{AChild1}{AChild2}{1}%
      \angled{BChild1}{BChild2}{1}%
    }, \\
    i'(\gendend)\succ i'(\gendend) &=
    \betaiw(\generator)\assocprod_\lambda \generator = \betagenerator
    \assocprod_\lambda \generator =     
    \tree{\children{A}{,xbbh=\lblheight,xbbl=.4\treesepone,xbbd=.7\levelsep}{2}{.8\treesepone}{.7\levelsep}%
      \rput[t](AChild1){\children{B}{}{2}{.8\treesepone}{.7\levelsep}}%
      \lbl[bl]{ARoot}{0pt}{\lblsep}{0}%
      \angled{AChild1}{AChild2}{1}%
      \angled{BChild1}{BChild2}{1}%
    }.
  \end{align*}
  Similarly we can prove that $i'(x\prec x) = i'(x)\prec i'(x)$ and
  $i'(x\ast x) = i'(x)\ast i'(x)$ in the case~$x=\gendend$.

  The inductive case is similar. Given two trees $x$, $y$ in
  $\FDT^\lambda$, we have $i'(x\succ y) = i'\bigl(\Gr(x\ast y_1,
  y_2,\ldots, y_m) \bigr)$, where $\H(y)=(y_1,\ldots,y_m)$. Note that
  the function $i'$ commutes with the grafting operations, in the sense
  that
  \begin{equation*}
    i'\bigl(\Gr(t_1,\ldots,t_k) \bigr) = \Gr\bigl( \betaiw i'(t_1),\ldots,
    \betaiw i'(t_k); 1,\ldots, 1 \bigr),
  \end{equation*}
  assuming that $i'(\point)=\point$, since the function $(f^0)^{-1}$
  and the normalization of the grafting produce the same results when
  collapsing the intermediate leaves of the trees.  Therefore, we can write
  \begin{align} 
    i'(x\succ y) &= i'\bigl(\Gr( x\ast y_1,\ldots y_m) \bigr) = 
    \Gr\bigl(\betaiw i'(x\ast y_1), \ldots, \betaiw i'(y_m); 1,\ldots,1\bigr) \notag \\
    &=\Gr\Bigl(\betaiw\bigl(i'(x)\ast i'(y_1)\bigr), \ldots, \betaiw i'(y_k);
    j_1,\ldots,j_{k-1} \Bigr) \label{E:xsuccy}
  \end{align}
  where in the last equality we have collapsed the intermediate leaves
  that are children of the root of $y$. Here we used the inductive
  hypothesis that guarantee $i'(x\ast y_1)=i'(x)\ast i'(y_1)$, as
  $y_1$ has less nodes than $y$.  On the other hand, to compute
  $i'(x)\succ i'(y) = \betaiw\bigl(i'(x)\bigr)\assocprod_\lambda i'(y)$,
  observe that $\overline{\betaiw i'(x)} = i'(x)$, since
  $i'(x)$ has root label $0$. Also, the de-grafting
  $\H\bigl(i'(y)\bigr)$ yields the same subtrees that appear
  in~\eqref{E:xsuccy}. Thus, using the definition of $\assocprod_\lambda$
  in~\eqref{E:def-prod}, we conclude that $i'(x)\succ i'(y)$ coincides
  with~\eqref{E:xsuccy}.  The other operations can be verified
  similarly.

  The map $i'$ sends $\gendend$ to $\generator$, thus $i'=i$ as
  claimed. The injectivity of $i'$ follows easily from the fact that
  it maps a linear basis of $\FDT^\lambda$ onto a subset of the linear
  basis $\tiw^0$ of~$\Fiw^\lambda$.
\end{proof}

According to Proposition~\ref{P:dendri-baxter}, a \rba\ of weight $0$ may be turned into a $0$-dendriform trialgebra, which as explained in Section~\ref{S:dend-free} gives rise to a dendriform dialgebra. Therefore, there is a unique
morphism of dendriform dialgebras
\begin{equation*}
\FDD \to \Fii^0
\end{equation*}
that sends $\gendend$ to $\generator$.
It is known that this map is injective~\cite{EG}. In fact, we can show
 that the composite
 \begin{equation*}
j:\FDD \to \Fii^0\onto \Fww^0
\end{equation*}
is still injective.

\begin{proposition}\label{P:dialg-inj} 
The canonical morphism of dendriform dialgebras
    \begin{equation*}
      j:\FDD\to\Fww^0,
    \end{equation*}
  sends any rooted planar binary tree $x$ to itself with root label $0$.
 In particular, $j$ is injective.
\end{proposition}
\begin{proof} The map $j$ is the composite of the following canonical maps:
  \begin{equation*}
\FDD\inc\FDT^0\map{i}\Fiw^0\onto\Fww^0\,.
\end{equation*}
The first map in this chain is the unique morphism of dendriform dialgebras
preserving the element $\gendend$. It is easy to see from the description of the operations in $\FDD$ and $\FDT^0$ that this map is simply the linearization of the inclusion of the set of rooted planar binary trees in the set of rooted planar trees. Applied to binary trees, the map $i$ merely adds a label $0$ to the root and a label $1$ to each angle of the tree, but does not change the underlying tree. The last map in the chain was described in Proposition~\ref{P:can-maps}; on trees with root label $0$ it simply erases the angle labels. Therefore $j$ is as claimed.
\end{proof}

\begin{remark}
By Proposition~\ref{P:trialg-inj},
 the image of the map $i:\FDT^\lambda\inc\Fiw^\lambda$ is precisely the subspace  $\Fiw^{\lambda,0}$ of  $\Fiw^{\lambda}$ spanned by $\tiw^0$. According to Proposition~\ref{P:filtration},  $\Fiw^{\lambda,0}$ is an ideal for the product $\assocprodlambda$ and a subalgebra for the
  product $\astlambda$ of $\Fiw^{\lambda}$. Moreover, that proposition implies that $\Fiw^{\lambda,0}$ is closed under the dendriform operations of $\Fiw^{\lambda}$. Thus, $i$ identifies $\FDT$ with the dendriform subtrialgebra  $\Fiw^{\lambda,0}$ of $\Fiw^{\lambda}$. This describes the free dendriform trialgebra explicitly as a subobject of the free \rba.
  
  The map $j:\FDD\inc\Fww^0$ of Proposition~\ref{P:dialg-inj} embeds
  the free dendriform dialgebra in the dendriform subdialgebra $\Fww^{0,0}$
  of $\Fww^{\lambda}$, but its image is strictly smaller.
  %  This
%  explains that the coefficient in the quasi-idempotent map
%  $\betaiw$~\eqref{E:convention-j} does not affect the result, and it
%  holds for arbitrary $\lambda$.
\end{remark}

\subsection{Dendriform dimensions v.s.\@ \rb\ dimensions}\label{S:vs}

Let $\FDT(n,m)$ be the the subspace of the free $\lambda$-dendriform trialgebra $\FDT^\lambda$ spanned by the set $\PT(n,m)$ (Section~\ref{S:comb-trees}). In other words, a rooted planar tree $x$ has  $\deg(x)=(n,m)$
if it has $n+1$ leaves  and $m$ internal nodes. In particular, $\deg(\gendend)=(1,1)$. 

This defines a bigrading on $\FDT^\lambda$  with similar properties to those of the bigrading of the free \rba\ $\Fiw^\lambda$ (Section~\ref{S:comb-dim}). Namely, the dendriform operations preserve the grading defined by the subspaces $\FDT(n,*)$ and the filtration defined by the subspaces
$\bigoplus_{\ell\leq m}\FDT(*,\ell)$. The morphism $i:\FDT^\lambda\to\Fiw^\lambda$ preserves the former grading and decreases the latter filtration degree by $1$,
 since according to Propositions~\ref{P:bijections} and~\ref{P:trialg-inj}, $i$ sends $\PT(n,m+1)$ to $\T_{\infty,2}(n,m)$ for $n\geq 1$, $m\geq 0$.

Consider the dimensions
\begin{gather*}
\dt(n,m)=\dim_{\field}\FDT(n,m),\quad \dt(n,*)=\dim_{\field}\FDT(n,*),
\\
\text{and}\quad \dt(k)=\dim_{\field}\bigoplus_{\substack{n,m\geq 1 \\n+m=k
  }}\FDT(n,m).
\end{gather*}
In view of Corollary~\ref{C:bijections}, we have the following relation between the  dimensions of the homogeneous components of
$\FDT^\lambda$ and $\Fiw^\lambda$:
\begin{gather*}
\fiw(n,m)=\dt(n,m)+\dt(n,m+1),\quad \fiw(n,*)=2\dt(n,*) \\
\text{and}\quad
\fiw(k)=\dt(k)+\dt(k+1)\,.
\end{gather*}
The first of these relations can be used to deduce the somewhat
complicated expression for $\dt(n,m)$~\eqref{E:plan-trees} from the simpler
expression for $\fiw(n,m)$ (Proposition~\ref{P:fiw}). The second one expresses the relation between the small and the large Schr\"oder numbers, while the last one
relates $\dt(k)$ to the Motzkin numbers (Proposition~\ref{P:fiw-k}).

%\fbox{Is $\dt(k)$ something else combinatorially?}

\smallskip

\enlargethispage{10pt}
We compare the dimensions of the free dendriform dialgebra $\FDD$ to those of the free \rba\ $\Fww^0$.
For a planar binary tree $x$  we have a notion of degree,
namely, $\deg(x)=n$ if $x$ has $n+1$ leaves. Such a tree has
$n$ angles and $n$ internal nodes, so $j(x)\in\tww(n,n-1)$ (the root of $j(x)$ has label $0$). It is well-known that the number of such trees is
the  Catalan number $C(n)$.
Since the map $j$ is injective, it follows from Proposition~\ref{P:fww} that
\begin{equation*}
C(n) \le nC(n-1),
\end{equation*}
for $n\ge 1$. We may view the embedding $j:\FDD\inc\Fww^0$ as an
algebraic realization of this inequality.  

\newpage 
\appendix
\section{Algorithms} \label{S:appendix}

The following algorithms are used in the proof of
Proposition~\ref{P:bijections}.

  \begin{algorithm}[!ht]
    \caption{$\treetopath{t}$: convert a tree $t\in\PT(n,m)$ to a path
      $p\in\SPp(n,m)$.}
    \label{treetopath}
    \begin{algorithmic}
    \IF{$t$ is a leaf}
      \STATE do nothing and return
    \ENDIF
    \STATE $\{t_1,\ldots,t_k\} \leftarrow \text{subtrees rooted at the
      children of the root
      of $t$, left to right}$
    \STATE write {\normalfont\sffamily H}
    \STATE $\treetopath{t_1}$
    \IF{$k>2$}
      \FOR{$i=2,\ldots,k-1$}
         \STATE write {\normalfont\sffamily D}
         \STATE $\treetopath{t_i}$
      \ENDFOR
    \ENDIF
    \STATE write {\normalfont\sffamily V}
    \STATE $\treetopath{t_k}$
  \end{algorithmic}
\end{algorithm}

\newcommand\pathtotree[1]{\operatorname{\mathsf{PathToTree}}(#1)}
\begin{algorithm}[!ht]
  \caption{$\pathtotree{s_1s_2\cdots s_k}$: convert a path
    $p\in\SPp(n,m)$ to a tree $t\in\PT(n,m)$.}
  \label{pathtotree}
  \begin{algorithmic}
%    \IF{$s_k$ is \textsf{D}}
%      \STATE return $\pathtotree{\mathsf{H}s_2\cdots
%        s_{k-1}\textsf{V}}$, with root label $0$
%    \ELSIF{the path has a diagonal step lying on the line $y=x$}
%      \STATE $s_i \leftarrow\text{first diagonal step over $y=x$}$
%      \STATE $s_j \leftarrow\text{first vertical step after $s_i$
%        which touch the line $y=x$}$
%      \STATE return $\pathtotree{\textsf{H}s_1\cdots
%        s_{i-1}\textsf{V}s_{i+2}\cdots s_{j-1}\textsf{H}s_{j+1}\cdots s_k\textsf{V}}$, with label $0$
%    \ELSE
      \STATE $t\leftarrow \text{\sffamily root}$
      \STATE $\textsf{node} \leftarrow \text{root of $t$}$
      \FOR{$i=1,\ldots,k$}
        \IF{$s_i$ is \textsf{H}}
          \STATE create a child \textsf{c} of \textsf{node} and mark
          \textsf{node} as \textsf{available}
        \ELSE
          \STATE $\textsf{node} \leftarrow \text{first parent of
          \textsf{node} with label \textsf{available}}$
          \STATE create a (rightmost) child \textsf{c} of \textsf{node}
          \IF{$s_i$ is \textsf{D}}
            \STATE mark \textsf{node} as \textsf{available}
          \ELSE
            \STATE mark \textsf{node} as \textsf{not available}
          \ENDIF
        \ENDIF
        \STATE $\textsf{node} \leftarrow \textsf{c}$
      \ENDFOR
      \STATE return the tree $t$
%   \ENDIF
  \end{algorithmic}
\end{algorithm}

%The following table lists the transformations used in the proof of part (ii) of Proposition~\ref{P:fww}, to relate the paths of part (i) to those of part (ii)
%(colored Motzkin paths).

\nocite{poisson,
B,
BSS,
Ca,
E,
EG,
GK,
GS,
L,
LR,
R69,
R1,
R2,
Sl,
St99,
catadd}

\newpage 
\bibliographystyle{plain}
\bibliography{db}

\end{document}